\documentclass[preprint,12pt]{article}
\usepackage{amsmath, amssymb, latexsym, amscd, amsthm,amsfonts,amstext}
\usepackage[mathscr]{eucal}
\usepackage{graphicx}
\usepackage{subfig}
\usepackage{float}
\usepackage{color}
\usepackage{hyperref}
\usepackage[utf8]{inputenc}
\usepackage[english]{babel}
\usepackage{algorithm}
\usepackage{algpseudocode}
\newtheorem{theorem}{Theorem}[section]

\newtheorem{problem}{Problem}[section]

\newtheorem*{assumption*}{Assumption}

\newtheorem{corollary}{Corollary}[section]
\newtheorem{remark}{Remark}[section]
\numberwithin{equation}{section}
\newcommand{\uinc}{u_{\rm inc}}

\newcommand{\ds}{\displaystyle}

\newcommand{\Pm}{P_{\rm meas}}
\newcommand{\Pp}{P_p}

\newcommand{\usc}{u_{\rm sc}}

\newcommand{\ii}{{\rm i}}
\newcommand{\ik}{\ii k}

\newcommand{\x}{{\bf x}}

\input{tcilatex}

\begin{document}

\title{A numerical method to solve a phaseless coefficient inverse problem from a single measurement of experimental data.}

\author{Michael V.Klibanov\thanks{
	Department of Mathematics and Statistics, 
	University of North Carolina at Charlotte, 
	Charlotte, NC, 28223, 
	nkoshev@uncc.edu, 
	mklibanv@uncc.edu, 
	lnguye50@uncc.edu.
} \and Nikolay A. Koshev\footnotemark[1] 
\and Dinh-Liem Nguyen\thanks{
	Department of Mathematics, 
	Kansas State University, Manhattan, KS 66506,
	dlnguyen@ksu.edu.
} \and Loc H. Nguyen\footnotemark[1] 
\and Aaron Brettin\thanks{
	Department of Physics and Optical Science,
	University of North Carolina at Charlotte, 
	Charlotte, NC, 28223, 
	abrettin@uncc.edu,
	astratov@uncc.edu.
}
\and Vasily N. Astratov\footnotemark[3]
}

%
%
%
\date{}
\maketitle
\begin{abstract}
	We propose in this paper a globally numerical method to solve a phaseless coefficient inverse
problem: how to reconstruct the spatially distributed refractive index of scatterers from the
intensity (modulus square) of the full complex valued wave field at an array of light detectors
located on a measurement board. The propagation of the wave field is governed by the 3D
Helmholtz equation. Our method consists of two stages. On the first stage, we use asymptotic analysis to obtain an upper estimate for
the modulus of the scattered wave field. This estimate allows us to approximately reconstruct
the wave field at the measurement board using an inversion formula. This reduces the
phaseless inverse scattering problem to the phased one. At the second stage, we apply a
recently developed globally convergent numerical method to reconstruct the desired refractive
index from the total wave obtained at the first stage. Unlike the optimization approach, the
two-stage method described above is global in the sense that it does not require a good initial
guess of the true solution. We test our numerical method on both computationally simulated
and experimental data. Although experimental data are noisy, our method produces quite
accurate numerical results. 
\end{abstract}




\section{Introduction}

\label{sec:1}

Using the apparatus of the Riemannian geometry and asymptotic analysis, we
construct in this paper a new numerical method for the solution of a 3D
Phaseless Coefficient Inverse Problem (phaseless CIP) for the Helmholtz
equation with the data resulting from a single measurement event. The
unknown coefficient of the Helmholtz equation is $n^{2}\left( \mathbf{x}%
\right) ,\mathbf{x}\in \mathbb{R}^{3},$ where $n\left( \mathbf{x}\right)
\geq 1$ is the spatially distributed refractive index. Our method computes
locations and refractive indices of unknown scatterers using experimentally
measured intensity (i.e. the square modulus) of the complex valued wave
field. That wave field is the solution of the Helmholtz equation. The phase
was not measured. We verify the accuracy of our computations via applying
the same numerical method to computationally simulated data.

Measurements were conducted on a part of a plane outside of scatterers. Only
a single direction of the incident plane wave on many wavelengths was used,
which means a single measurement event. This is more difficult than the case
of multiple measurements. The experimental data were collected by ourselves.
The authors are unaware about other publications in which a 3D phaseless CIP would
be computationally solved for the case when the experimental intensity data
would be collected on several wavelengths and for a single measurement event.

Our scatterers are microspheres of the diameter of 6 $\mu m$ (micron). Our
experimental data were collected for the case of the vertically propagated
white light. To obtain data on an interval of frequencies, the light was
filtered on six (6) wavelengths ranging from 0.420$\mu m$ to 0.671$\mu m.$
Following, e.g., \cite{Tseng}, to measure the light intensity, we have used
the detector array which is available in the camera of the Samsung Galaxy S3 mobile
telephone unit. The idea of \cite{Tseng} is to built an extremely low weight
optical system.

We used white light source with a set of narrow band filters to provide various 
illumination wavelengths because this technique can be used in combination with 
standard microscopes and in any environment such for example as clinical environment. This way we have measured the
intensity of only the full wave field. However, since the case of the full
wave field was not considered in the above cited publications on phaseless CIPs, we
develop here a significantly new numerical method for our phaseless CIP. 
Phaseless CIPs for the case when the intensity of the full wave field is measured
were also studied analytically in \cite{RomYam} for both Helmholtz and
Schr\"odinger equations. In this case the medium is simultaneously illuminated
by two point sources and many such pairs of sources are used. This
does not work, however, for our experimental arrangement, since we have only a single
direction of the incident plane wave.

First, we establish an  
inversion formula to  approximate
 the wave field at the measurement site. This is considered as the first stage of our
numerical method. 
 This inversion formula is very interesting because the reconstruction of the complex number $z = |z|e^{\ii \arg z}$ from $z$ is not unique. In this paper, we can successfully derive this approximate reconstruction by using an asymptotic behavior of the wave field and proving {\it a priori} bound for the scattering wave.
After this stage, we obtain a Phased Coefficient Inverse Problem (phased CIP). 
 On the second stage, we solve that phased CIP to
reconstruct the unknown coefficient of the Helmholtz equation. It is on the
second stage when we reconstruct locations and refractive indices of those
microspheres. The numerical solution of the phased CIP is found using the
globally convergent numerical method, which was recently developed in \cite%
{KlibanovLiemLocHui:IPI2017}, also, see \cite%
{AlekKlibanovLocLiemThanh:anm2017,LiemKlibanovLocAlekFiddyHui:jcp2017,nguyen2017:iip2017}
for the performance of this method on microwave experimental backscattering
data

Our interest in phaseless CIPs is motivated by applications to optical imaging of
such small objects as, e.g. biological cells and microspheres. Sizes of
biological cells are usually in the interval of $\left( 5,100\right) \mu m$ 
\cite{Phill2009}. To optically image such small targets, one should use
light sources, in which case the wavelengths are of $1\mu m$ range. The
wavelength $\lambda =1\mu m$ corresponds to the frequency $\omega \approx
300,000$ Gigahertz. It is currently impossible, however, to arrange stable
measurements of the phase for such high frequencies. Only the intensity of
the scattered wave field can be reliably measured on these frequencies \cite%
{Petersena:U2008, Ruhlandt:PR2014}.

While we have measured the intensity of the full wave field, it is also
possible sometimes to measure the intensity of the scattered wave field. On
the other hand, a number of past works for phaseless CIPs of the first author with
coauthors were devoted to the analytical reconstruction procedures for the
case when the intensity of the scattered wave field is measured \cite%
{KR:JIpp2015,KlibanovRomanov:ejmca2015,KlibanovRomanov:SIAMam2016,KlibanovRomanov:ip2016,KlibanovLiemLoc:arxiv2017}%
. To arrange experimental measurements for this case, one needs to work with
tunable lasers, which would operate on several wavelengths.

The question on how to solve the inverse scattering problem without the
phase information was probably first posed in the book of Chadan and
Sabatier \cite[Chapter 10]{Chadan:sv1977} published in 1977. Fifteen years
later, the first uniqueness result for this problem in the 1D case was
established in \cite{Klibanov:jmp1992}, also see \cite{AktosunSacks:ip1998}
for a follow up result. Next, the first uniqueness result in 3D was obtained
in \cite{Klibanov:SiamJMath2014}. Since then, the 3D phaseless CIPs were studied
intensively.
In \cite{KlibanovLocKejia:apnum2016} a modified reconstruction procedure of 
\cite{KlibanovRomanov:SIAMam2016} was numerically implemented. Multiple
locations of the point source at multiple frequencies were used in \cite%
{KlibanovRomanov:SIAMam2016,KlibanovLocKejia:apnum2016}. Unlike this, in 
\cite{KlibanovLiemLoc:arxiv2017} the case of a single direction of the
incident plane wave on a frequency interval was numerically implemented. 
 We refer the reader to other versions of the uniqueness
theorems for 3D phaseless CIPs in \cite%
{Klibanov:AppliedMathLetters2014,Klibanov:aa2014,Klibanov:ipi2017,KlibanovRomanov:ip2017}%
. As mentioned above, the analytic reconstruction procedures in for 3D
phaseless CIPs were proposed in \cite%
{KR:JIpp2015,KlibanovRomanov:ejmca2015,KlibanovRomanov:SIAMam2016,KlibanovRomanov:ip2016}
for the \ case when the intensity of the scattered rather than full wave
field is measured.

We now refer to other approaches to phaseless inverse scattering problems.
In \cite{BardsleyVasquez:sjis2016, BardsleyVasquez:ip2016} a phaseless CIP for
Helmholtz equation was solved numerically using Kirchhoff migration and Born
approximation. While coefficients of partial differential equations are
subjects of interests in all above cited works, there is also a significant
interest in the reconstruction of surfaces of scatterers from the phaseless
data. In this regard we refer to publications \cite%
{AmmariChowZou:sjap2016,Bao2013,BaoZhang:ip2016,IvanyshynKressSerranho:acm2010,Ivany2010,Ivany2011, LiLiuWang:jcp2014}
and references cited therein. We also mention the problem of the
reconstruction of a compactly supported function from the absolute value of
its Fourier transform \cite{Klib2006,KlibKam} as well as a closely related
problem of the solution of the autoconvolution equation \cite%
{Gerth2014,Hofmann2,Burge2016}.

As to the phased CIPs, they arise in many real world applications including
detection and identification of explosives, non-destructive testing, medical
imaging, and geophysics prospecting. In general, CIPs are nonlinear and
ill-posed. The developments of
the numerical methods for CIPs are challenging. Due to a large variety of applications,
there is a huge literature on numerical reconstruction methods for these
problems. \ We refer here to a few publications and references cited therein 
\cite%
{Ammar2014,Ammar2013,AmmariChowZou:sjap2016,BeilinaKlibanovBook,ChowZou,ColtonKress:2013,Gonch2013,Gonch2017,Lakha2010,Ito2013, LiLiuWang:jcp2014, Li2015, Li2017}.

We call a numerical method for a CIP \emph{globally convergent} if there is
a theorem, which claims that this method delivers at least one point in a
sufficiently small neighborhood of the exact solution without any advanced
knowledge of this neighborhood. As to the numerical method of \cite%
{KlibanovLiemLocHui:IPI2017}, which is used here, Theorem 6.1 of \cite%
{KlibanovLiemLocHui:IPI2017} ensures its global convergence. While the above
mentioned globally convergent numerical method works for CIPs with single
measurement data, we also refer to \cite{Kaban2004,Kaban2015,Kaban2015a} for
a global reconstruction technique for a CIP with the data resulting from
multiple measurements. The idea of these references is based on an extension
to the 2D case of the well-known Gelfand-Krein-Levitan method, which works
for a 1D CIP.

In Section \ref{sec:2}, we state our phaseless CIP. In Section \ref{sec:3}, we
recall the asymptotic behavior of the solution of the Helmholtz equation
when the wave number tends to infinity. In Section \ref{sec:4} we
estimate both analytically and numerically the intensity of the scattered
wave field and, using this estimate derive inversion formulae which enable
us to approximate the wave field at the measurement site. In Section \ref%
{sec:5}, we briefly outline the globally convergent numerical method of \cite%
{KlibanovLiemLocHui:IPI2017}. In Section \ref{sec:6}, we describe the
procedure of the collection of the experimental data. In Section \ref{sec:7}
we present our numerical results. In Section \ref{sec:8} We summarize results 
of this paper in Section \ref{sec:8}.

\section{Problem statement}

\label{sec:2}

%

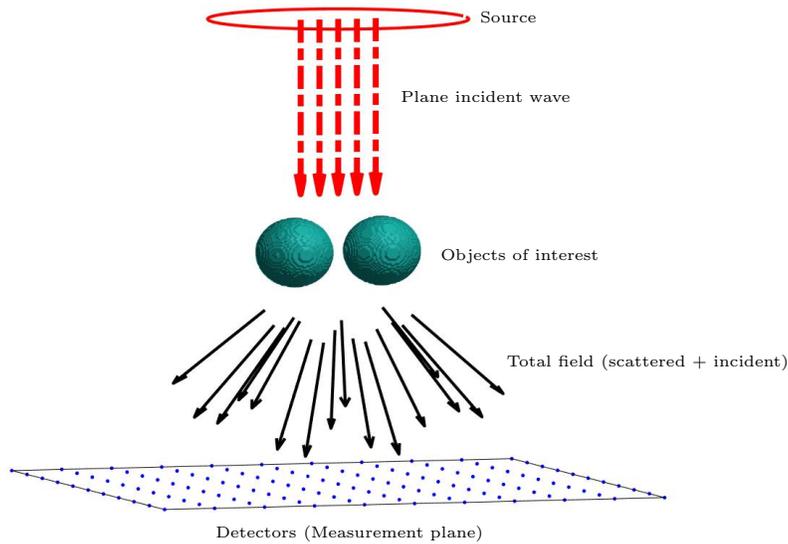
\begin{figure}[tbp]
	\begin{center}
		\includegraphics[width = 0.7\linewidth]{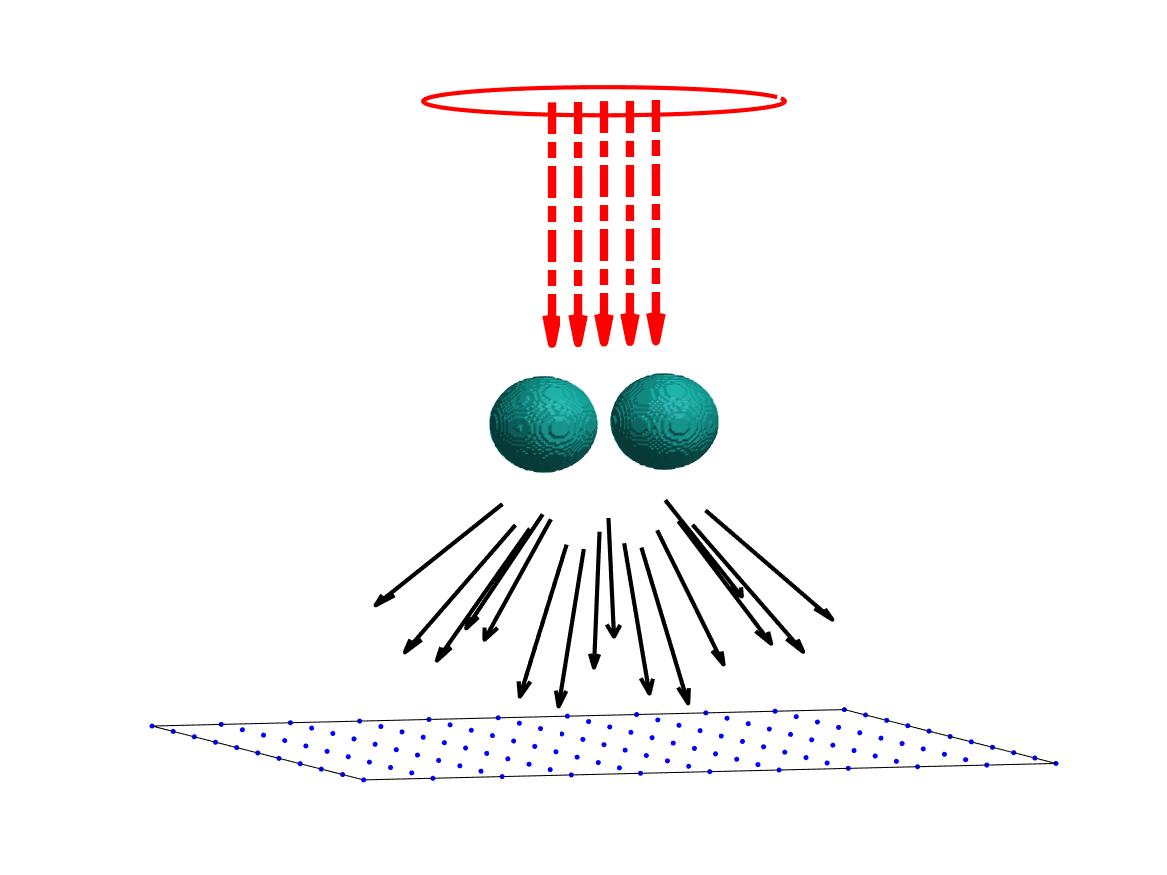}
		\put(-200,15){\tiny Detectors (Measurement plane)}
		\put(-90,80){\tiny Total field (scattered + incident)}
		\put(-115,120){\tiny Objects of interest}
		\put(-130,180){\tiny Plane incident wave}
		\put(-100,210){\tiny Source}
	\end{center}
		\caption{ \emph{A schematic diagram of measurements. }}
	\label{fig diagram}
\end{figure}

In this section we formulate the phaseless CIP of this paper. To this end we first
briefly describe the direct scattering problem. Suppose that an object is
illuminated by an incident plane wave. The interaction of this incident wave
with the object produces the scattered wave, see Figure \ref{fig diagram}.
The total wave field is the sum of the incident wave and the scattered wave.

We denote $\mathbf{x}=\left( x_{1},x_{2},x_{3}\right) $ points of $\mathbb{R}%
^{3}$. Let $\Omega \subset \mathbb{R}^{3}$ be a bounded domain with a
piecewise smooth boundary $\partial \Omega $, containing the scatterers. It
is convenient for our computational purpose to specify the domain $\Omega $
as%
\begin{equation}
\Omega =\left\{ \mathbf{x}:\left\vert x_{1}\right\vert ,\left\vert
x_{2}\right\vert <b,x_{3}\in \left( -d_{1},d_{2}\right) \right\},  \label{1}
\end{equation}%
where $b,d_{1},d_{2}>0$. We also denote $\Gamma $ one of
sides of $\Omega ,$%
\begin{equation}
\Gamma =\left\{ \mathbf{x}:\left\vert x_{1}\right\vert ,\left\vert
x_{2}\right\vert <b,x_{3}=d_{2}\right\} \subset \partial \Omega .  \label{2}
\end{equation}%
Let the function $n\left( \mathbf{x}\right) $, defined for all $\mathbf{x}%
\in \mathbb{R}^{3},$ represent the spatially distributed refractive index of
the medium. We assume that microspheres of our interest are located in the
domain $\Omega .$ Assume that 
\begin{equation}
n(\mathbf{x})=1\mbox{ for all }\mathbf{x}\in \mathbb{R}^{3}\setminus 
\overline{\Omega }\quad \mbox{ and }\quad n(\mathbf{x})\geq 1\mbox{ for all }%
\mathbf{x}\in \mathbb{R}^{3}.  \label{cond c}
\end{equation}%
Condition \eqref{cond c} means that the dielectric constant of the
background is scaled to be 1 and that of the scattering object is greater
than 1. Let $k>0$ be the wave number, consider the incident plane wave 
\begin{equation}
\uinc(\mathbf{x},k)=e^{\mathrm{i}kx_{3}}.  \label{20}
\end{equation}%
Denote by $u_{\mathrm{sc}}(\mathbf{x},k)$ the scattering wave. Then, the
 total wave field 
\begin{equation}
u(\mathbf{x},k)=\uinc(\mathbf{x},k)+u_{\mathrm{sc}}(\mathbf{x},k)
\label{eqn usc}
\end{equation}%
is governed by the Helmholtz equation with the Sommerfeld outgoing radiation  condition at the
infinity, 
\begin{equation}
\left\{ 
\begin{array}{rcll}
\Delta u(\mathbf{x},k)+k^{2}n^{2}(\mathbf{x})u(\mathbf{x},k) & = & 0 & 
\mathbf{x}\in \mathbb{R}^{3}, \\ 
\partial _{|\mathbf{x}|}u_{\mathrm{sc}}(\mathbf{x},k)-\mathrm{i}ku_{\mathrm{%
sc}}(\mathbf{x},k) & = & o(|\mathbf{x}|^{-1}) & |\mathbf{x}|\rightarrow
\infty .%
\end{array}%
\right.  \label{eqn Helmholtz}
\end{equation}%
Let the number%
\begin{equation}
R>d_{2}.  \label{21}
\end{equation}%
Define the plane $P$ and a square $P_{\text{meas}}\subset P$ as 
\begin{equation}
P=\{\mathbf{x=}(x_{1},x_{2},R):x_{1},x_{2}\in \mathbb{R}\}, \quad P_{\text{meas}}=\{%
\mathbf{x=}(x_{1},x_{2},R):\left\vert x_{1}\right\vert ,\left\vert
x_{2}\right\vert <b\}.  \label{3}
\end{equation}%
We call $P$ the \textquotedblleft measurement plane." Measurements of the
intensity are conducted on the square $P_{\text{meas}}$ 
for the wave numbers $k\in \lbrack \underline{k},\overline{k}].$ Here, the
interval $k\in \lbrack \underline{k},\overline{k}]$ represents the allowable
range of wave numbers.

\begin{problem}[The phaseless coefficient inverse scattering problem]
\label{Problem phaseless} Let $1\ll \underline{k}<\overline{k}<\infty $.
Assume that the function $f(\mathbf{x},k),$ 
\begin{equation}
f(\mathbf{x},k)=|u(\mathbf{x},k)|^{2},\quad \mathbf{x}\in P_{\text{meas}%
},k\in \lbrack \underline{k},\overline{k}]  \label{4}
\end{equation}%
is known. Determine the function $n(\mathbf{x})$ for $\mathbf{x}\in \Omega $.
\end{problem}

\begin{remark}[A comment on the Helmholtz equation]
Although the full Maxwell's system is the right model to describe the
propagation of the total wave field, we use the Helmholtz equation %
\eqref{eqn Helmholtz} in this paper. We have numerically verified in \cite[%
Section 8]{KlibanovLiemLoc:arxiv2017} that if the incident wave field has
the form $(0,u_{inc}(\mathbf{x},k),0)$, then $E_{2}(\mathbf{x},k)$, the
second component of the electric wave field satisfying the Maxwell's system,
matches well the total wave field $u(\mathbf{x},k)$. The study for the
phaseless coefficient inverse problem for Maxwell's system with general
incident wave field is considered as future research.
\end{remark}

It is well-known that the Helmholtz equation \eqref{eqn Helmholtz} can be
reformulated as the Lippmann-Schwinger equation (see \cite[Chapter 8]%
{ColtonKress:2013})

\begin{equation}
u(\mathbf{x},k)=e^{\text{i}kx_{3}}+k^{2}\displaystyle \int\limits_{\Omega }%
\frac{\exp (\mathrm{i}k|\mathbf{x}-\mathbf{\xi }|)}{4\pi |\mathbf{x}-\mathbf{%
\xi }|}(n^{2}(\mathbf{\xi })-1)u(\mathbf{\xi },k)d\mathbf{\xi ,}\quad 
\mathbf{x}\in \mathbb{R}^{3}.  \label{eqn LS}
\end{equation}%

Using the method in \cite{LechleiterNguyen:acm2014,Nguyen:anm2015}, we solve numerically the integral equation \eqref{eqn LS} to computationally
simulate the data for the phaseless CIP. In addition, we solve this equation
iteratively on the second stage of our reconstruction procedure.

\section{Asymptotic behavior of the total wave  as $k\rightarrow \infty $}

\label{sec:3}

In this section we establish the asymptotic behavior of the function $%
u\left( \mathbf{x},k\right) $ and at $k\rightarrow \infty .$ Although results of
this section follow from \cite{KlibanovRomanov:ip2016}, we need to formulate
these results again here since we essentially use them in our numerical
method.

\subsection{Geodesic lines}

\label{sec:3.1}

In addition to conditions (\ref{cond c}) imposed on the function $n\left( 
\mathbf{x}\right) ,$ we also assume that 
\begin{equation}
n\in C^{15}\left( \mathbb{R}^{3}\right) .  \label{7}
\end{equation}%
The smoothness condition (\ref{7}) is a technical one. It was used in \cite%
{KlibanovRomanov:ip2016} to prove an analog of Theorem 3.1. And that analog,
in turn was derived in \cite{KlibanovRomanov:ip2016} using the construction
of the solution of the Cauchy problem for a certain hyperbolic equation in 
\cite{Roman2002}. This construction technically needs (\ref{7}). In
addition, usually extra smoothness assumptions are not of a great concern in
the theory of CIPs, see, e.g. Theorem 4.1 in \cite{Roman1986}.

The Riemannian metric generated by the function $n\left( \mathbf{x}\right) $
is%
\begin{equation}
d\tau (\mathbf{x})=n(\mathbf{x})|d\mathbf{x}|,\quad |d\mathbf{x}|=\sqrt{%
(dx_{1})^{2}+(dx_{2})^{2}+(dx_{3})^{2}}.  \label{8}
\end{equation}%
Let the number $a>d_{1}.$ Consider the plane 
\begin{equation*}
P_{a}=\{\mathbf{x}=(x_{1},x_{2},-a):x_{1},x_{2}\in \mathbb{R}\}.
\end{equation*}%
Then by (\ref{1}) and (\ref{21}), both $\Omega$ and $P$ are contained in $\left\{ x_{3}>-a\right\}
.$ We assume below without further mentioning the condition about the regularity of geodesic
lines:

\begin{assumption*}[Regularity of geodesic lines]{\rm
\emph{For any point }$\mathbf{x}\in 
\mathbb{R}^{3}$\emph{\ there exists a unique geodesic line }$\Gamma (\mathbf{%
x},a)$\emph{, with respect to the metric }$d\tau $\emph{\ in (\ref{8})
connecting }$\mathbf{x}$\emph{\ with the plane }$P_{a}$\emph{\ and
perpendicular to }$P_{a}$\emph{.} }
\end{assumption*}
Again, we need this assumption only for Theorem \ref{thm 1}. But we do not verify it
in our numerical studies. A sufficient condition of the
regularity of geodesic lines is  (see \cite{Roman2014}) 
\begin{equation*}
\displaystyle \sum\limits_{i,j=1}^{3}\frac{\partial ^{2}\ln \left( n(\mathbf{%
x})\right) }{\partial x_{i}\partial x_{j}}\xi _{i}\xi _{j}\geq 0,\quad \text{%
for all }\mathbf{x}\in \overline{\Omega },\mathbf{\xi }\in \mathbb{R}^{3}.
\end{equation*}%
For an arbitrary point $\mathbf{x}\in \left\{ x_{3}>-a\right\} $ the travel
time along the geodesic line $\Gamma (\mathbf{x},a)$ from the plane $P_{a}$
to the point $\mathbf{x}$ is (see \cite{KlibanovRomanov:ip2016}) 
\begin{equation}
\tau (\mathbf{x})=\displaystyle \int\limits_{\Gamma (\mathbf{x},a)}n\left( 
\mathbf{\xi }\right) d\sigma .  \label{9}
\end{equation}%
The following theorem follows immediately from formulae (4.24)-(4.26) of 
\cite{KlibanovRomanov:ip2016}:

\begin{theorem}
Let $G\subset \left\{ x_{3}>-a\right\} $ be an arbitrary bounded domain such
that $\Omega ,P_{\text{meas}}\subset G.$ Then the following asymptotic
behavior of the solution of the problem (\ref{eqn usc}), (\ref{eqn Helmholtz}%
) holds 
\begin{equation}
u(\mathbf{x},k)=A(\mathbf{x})\exp (\mathrm{i}k\tau (\mathbf{x}))\left( 1+g(%
\mathbf{x},k)\right) ,\quad \mathbf{x}\in G,k\rightarrow \infty ,
\label{eqn asym}
\end{equation}%
where functions $A(\mathbf{x})>0$ and $g(\mathbf{x},k)$ are smooth and also
for $j=1,2,3$ 
\begin{equation}
g(\mathbf{x},k),\partial _{k}^j g(\mathbf{x},k)=O\left( 1/k\right)
,\quad \mathbf{x}\in G,k\rightarrow \infty .  \label{10}
\end{equation}%
\label{thm 1}
\end{theorem}

Corollary 3.1 follows immediately from \eqref{eqn asym}, \eqref{10}:

\begin{corollary}
We have 
\begin{equation}
A(\mathbf{x})=\lim_{k\rightarrow \infty }|u(\mathbf{x},k)|,\text{ }\forall 
\mathbf{x}\in G.  \label{110}
\end{equation}%
\label{corA}
\end{corollary}

Using (\ref{eqn asym}), we obtain the following approximate formula for the
scattered wave field $u_{\text{sc}}\left( \mathbf{x},k\right) $ for
sufficiently large values of the wave number $k$: 
\begin{equation}
u_{\text{sc}}\left( \mathbf{x},k\right) =A\left( \mathbf{x}\right) \exp
\left( \text{i}k\tau \left( \mathbf{x}\right) \right) -e^{\text{i}kx_{3}},%
\text{ }\mathbf{x}\in G.  \label{160}
\end{equation}

\section{An analytical upper estimate of $\left\vert u_{\text{sc}}\left( 
\mathbf{x},k\right) \right\vert ^{2}$ and an approximate inversion formula}

\label{sec:4}

\begin{figure}
	\centering
	\subfloat[\label{fig 2a}\emph{ Case of one microsphere described in Section \ref{sec:7.1}. In this case} $\protect \varphi 
		\left( k\right) \leq 0.038.$]{
		\includegraphics[width = 0.45\linewidth]{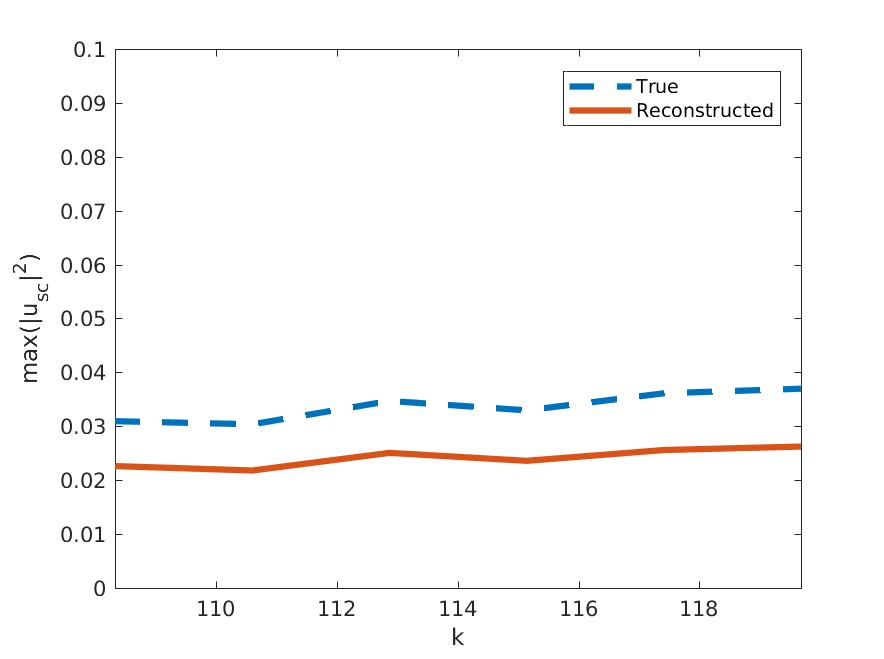}  
		}\hfill
	\subfloat[\label{fig 2b}\emph{ The case of two microspheres described in Section \ref{sec:7.2}. In this case $\protect\varphi \left(
		 k\right) \leq 0.115.$}]{
		\includegraphics[width = 0.45\linewidth]{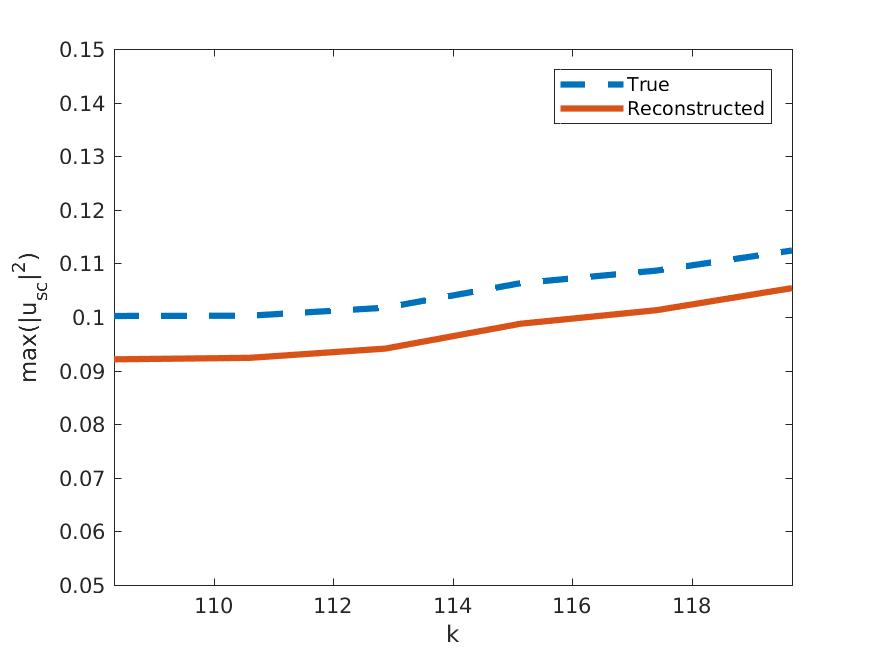} 
		}
	\caption{ \label{fig 2}\it The function $\protect\varphi \left( k\right)$ defined in \eqref{def varphi} is displayed for $k\in \left[ \protect%
	\underline{k}, \overline{k}\right] =\left[ 108.2,119.6\right]. $ It is computed from the
	 simulated data generated from exactly the same models and parameters as in the
	experiments. The
	function $\protect\varphi \left( k\right) $  is sufficiently small. Here, \textquotedblleft true" means the function $%
	\protect\varphi \left( k\right) $ results from the solution of
	equation (\protect\ref{eqn LS}). \textquotedblleft Reconstructed" means that the function is computed from 
	the reconstruction of the total wave $u(\x, k)$ using the method using formula \eqref{4.18} in this section 
	and \eqref{eqn usc} . }
\end{figure}

To derive our approximate inversion formula, we need to assume that the
function $\left\vert u_{\text{sc}}\left( \mathbf{x},k\right) \right\vert
^{2} $ is sufficiently small. We have computationally estimated $\left\vert
u_{\text{sc}}\left( \mathbf{x},k\right) \right\vert ^{2}$ for those
parameters which we use in our studies of both computationally simulated and
experimental data. Indeed, let $k$ be the dimensionless wave number, see
Section \ref{sec:6} for details. We have
measured data for $k\in \left[ 93.59,149.52\right] ,$\emph{\ }which
corresponds to the wavelengths $\lambda \in \left[ 0.420,0.671\right] \mu m$. However, we have observed that
 data are too noisy for $k\notin \left[ 108.3,119.6\right] ,$ see (\ref%
{6.2}). Hence, we set%
\begin{equation}
k\in \left[ \underline{k},\overline{k}\right] =\left[ 108.3,119.6\right] .
\label{4.1}
\end{equation}%
We have computationally simulated the data for the interval of wave numbers (%
\ref{4.1}) via the numerical solution of equation (\ref{eqn LS}). In doing
so, we have modeled the microspheres by exactly the same parameters as they
are in the experiment. Figures \ref{fig 2a}--\ref{fig 2b} display the graph of the function
\begin{equation}
	\varphi(k) = \max_{\x \in \Pm} |\usc(\x, k)|^2 \quad k \in [\underline k, \overline k].
	\label{def varphi}
\end{equation}
 These illustrations indicate that our assumption about the smallness of the function $\varphi(k)$ is might be true, at least for models and the range of the parameters used in this paper.

\subsection{An analytical upper estimate}

\label{sec:4.1}

It is desirable to provide an analytical justification to our computational
finding (Figure \ref{fig 2}) an upper bound for the function $\varphi(k),$ $k\in \left[ \underline{k},\overline{k}\right] .$ It follows immediately from (\ref{20}), (\ref{eqn usc}) and (\ref%
{eqn LS}) that
\begin{equation}
u_{\text{sc}}(\mathbf{x},k)=k^{2}\displaystyle \int\limits_{\Omega }\frac{%
\exp (\mathrm{i}k|\mathbf{x}-\mathbf{\xi }|)}{4\pi |\mathbf{x}-\mathbf{\xi }|%
}\left( n^{2}(\mathbf{\xi })-1\right) u(\mathbf{\xi },k)d\mathbf{\xi ,}\text{
}\mathbf{x}\in P_{\text{meas}}.  \label{4.2}
\end{equation}%
We can see from \eqref{4.2} that $|\usc(\x, k)|$ grows at the order $O(k^2)$ as $k \to \infty$. Our
goal in Theorem \ref{thm 4.1} is to decrease that growth to be $O(k)$. Since by (%
\ref{3}), $\mathbf{x=}\left( x_{1},x_{2},R\right) $ on $P_{\text{meas}}.$
Assuming that $R \gg \max \left( b,1\right) $ and using the well-known
formula for the far field approximation, we obtain from (\ref{4.2})%
\begin{equation}
u_{\text{sc}}(\mathbf{x},k)=\frac{k^{2}}{4\pi R}e^{\text{i}kR}\displaystyle %
\int\limits_{\Omega }e^{-\text{i}k\xi _{3}}\left( n^{2}(\mathbf{\xi }%
)-1\right) u(\mathbf{\xi },k)d\mathbf{\xi ,} \quad \mathbf{x}\in P_{\text{%
meas}}.  \label{4.3}
\end{equation}%
We use in (\ref{4.3}) \textquotedblleft $=$" instead of \textquotedblleft $%
\approx $" only for the convenience of the {\it a priori} bound established in Theorem \ref{thm 4.1} below.

\begin{theorem}
\rm \emph{Assume that} 
\begin{equation}
\left\vert \tau _{\xi _{3}}\left( \mathbf{\xi }\right) -1\right\vert \geq
\beta =const.>0, \quad \forall \mathbf{\xi }\in \Omega .  \label{4.4}
\end{equation}%
\emph{Then the following estimate is valid for the function} $u_{\text{sc}}(%
\mathbf{x},k)$ \emph{given in (\ref{4.3})}%
\begin{equation}
\left\vert u_{\text{sc}}(\mathbf{x},k)\right\vert \leq C\frac{k}{4\pi R}%
\left( 1+O\left(1/k\right) \right), \quad \text{ }\mathbf{x}\in P_{\text{%
meas}}, k\rightarrow \infty ,  \label{4.5}
\end{equation}%
\emph{where the number }$C>0$\emph{\ depends only on the domain }$\Omega ,$%
\emph{\ functions} $A\left( \mathbf{\xi }\right) ,c\left( \mathbf{\xi }%
\right) $ for $\mathbf{\xi }\in \Omega $ \emph{and} \emph{the number} $\beta 
$ \emph{in (\ref{4.4}). }
\label{thm 4.1}
\end{theorem}

\textbf{Proof}. Since the function $n^{2}(\mathbf{x})-1=0$ outside of the
domain $\Omega $, then by (\ref{1}), we can assume that $n^{2}(\mathbf{x}%
)-1=0$ for $\left\vert \mathbf{x}\right\vert >r,$ where the number $r>0$ is
such that $B\left( r\right) =\left\{ \left\vert \mathbf{x}\right\vert
<r\right\} \subset \Omega .$ Using $u(\mathbf{\xi },k)=A\left( \mathbf{\xi }%
\right) e^{\text{i}k\tau \left( {\bf \xi} \right) }\left( 1+O\left( 1/k\right)
\right) $ and substituting in (\ref{4.3}), we obtain
\begin{multline}
u_{\text{sc}}(\mathbf{x},k)=
\frac{k^{2}}{4\pi R}e^{\text{i}kR}\displaystyle\int\limits_{\xi _{1}^{2}+\xi
_{2}^{2}<r^{2}}d\xi _{1}d\xi _{2}
\displaystyle\int\limits_{-\sqrt{r^{2}-\xi
_{1}^{2}-\xi _{2}^{2}}}^{\sqrt{r^{2}-\xi _{1}^{2}-\xi _{2}^{2}}}e^{-\text{i}%
k\xi _{3}}e^{\text{i}k\tau \left( \mathbf{\xi }\right) }
\\
\left( 1+O\left(
1/k\right) \right) A\left( \mathbf{\xi }\right) \left( n^{2}(\mathbf{\xi }%
)-1\right) d\xi _{3}\mathbf{.}
\label{4.6}
\end{multline}
We now estimate the interior integral in (\ref{4.6}),%
\begin{equation}
I\left( \xi _{1},\xi _{2},k,r\right) =\displaystyle\int\limits_{-\sqrt{%
r^{2}-\xi _{1}^{2}-\xi _{2}^{2}}}^{\sqrt{r^{2}-\xi _{1}^{2}-\xi _{2}^{2}}%
}e^{-\text{i}k\xi _{3}}e^{\text{i}k\tau \left( \mathbf{\xi }\right) }\left(
1+O\left( 1/k\right) \right) A\left( \mathbf{\xi }\right) \left( c(\mathbf{%
\xi })-1\right) d\xi _{3}.  \label{4.7}
\end{equation}%
By (\ref{4.4}) we can assume without a loss of generality that 
\begin{equation}
\tau _{\xi _{3}}\left( \mathbf{\xi }\right) -1\geq \beta \text{ in }\Omega .
\label{4.8}
\end{equation}%
Change variables $\xi _{3}\Leftrightarrow z=\tau \left( \mathbf{\xi }\right)
-\xi _{3}.$ By the implicit function theorem and (\ref{4.8}) this equation
can be uniquely solved with respect to $\xi _{3}$ as $\xi _{3}=\xi
_{3}\left( z,\xi _{1},\xi _{2}\right) .$ Then 
\begin{equation}
d\xi _{3}=\frac{dz}{\tau _{\xi _{3}}\left( \xi _{1},\xi _{2},\xi _{3}\left(
z,\xi _{1},\xi _{2}\right) \right) -1}.  \label{4.9}
\end{equation}%
Denote $\rho =\rho \left( \xi _{1},\xi _{2},r\right) =\sqrt{r^{2}-\xi
_{1}^{2}-\xi _{2}^{2}}.$ Then (\ref{4.7})-(\ref{4.9}) imply that%
\begin{equation}
I\left( \xi _{1},\xi _{2},k,r\right) =\displaystyle\int\limits_{\tau \left(
\xi _{1},\xi _{2},-\rho \right) +\rho }^{\tau \left( \xi _{1},\xi _{2},\rho
\right) -\rho }e^{\text{i}kz}\left( 1+O\left( 1/k\right) \right) \frac{\left[
A\left( n^{2}-1\right) \right] \left( \xi _{1},\xi _{2},\xi _{3}\left( z,\xi
_{1},\xi _{2}\right) \right) }{\tau _{\xi _{3}}\left( \xi _{1},\xi _{2},\xi
_{3}\left( z,\xi _{1},\xi _{2}\right) \right) -1}dz.  \label{4.10}
\end{equation}%
Next, for any finite interval $\left( a,b\right) $ and for any complex
valued function $\phi \in C^{1}\left[ a,b\right] $%
\begin{equation*}
\displaystyle\int\limits_{a}^{b}\phi \left( x\right) e^{\text{i}kx}dx=\frac{1%
}{ik}\left( \phi \left( b\right) e^{\text{i}kb}-\phi \left( a\right) e^{%
\text{i}ka}\right) -\frac{1}{ik}\displaystyle\int\limits_{a}^{b}\phi
^{\prime }\left( x\right) e^{\text{i}kx}dx=O\left(1/k\right)
, \quad k\rightarrow \infty .
\end{equation*}%
Hence, using (\ref{10}), (\ref{4.8}) and (\ref{4.10}), we obtain%
\begin{equation}
I\left( \xi _{1},\xi _{2},k,r\right) =O\left( 1/k\right)
, \quad k\rightarrow \infty .  \label{4.11}
\end{equation}%
The inequality (\ref{4.5}) follows immediately from (\ref{4.6}), (\ref{4.9}) and (%
\ref{4.11}). \ $\square $

\subsection{Comparison with Figures \ref{fig 2}}

\label{sec:4.2}

Hence, by (\ref{4.5})%
\begin{equation}
\left\vert u_{\text{sc}}(\mathbf{x},k)\right\vert ^{2}\leq C\left( \frac{k}{%
4\pi R}\right) ^{2}\left( 1+O\left( 1/k\right) \right),\text{ }%
\mathbf{x}\in P_{\text{meas}}, \quad k\rightarrow \infty .  \label{4.12}
\end{equation}%
Consider now specific values of parameters $k$ and $R$ which we have used
both in computationally simulated and experimental data, substitute them in
the term $\left[ k/\left( 4\pi R\right) \right] ^{2}$ in (\ref{4.12}) and
compare resulting values in Figures \ref{fig 2}. Since we have had one and two
microspheres in the cases of Figures \ref{fig 2a} and \ref{fig 2b} respectively, we will
multiply that term by 2 in the case of two microspheres: $\left[ 2
k/\left( 4\pi R\right) \right] ^{2}$. Note that in our experiments and
computationally simulated data,  $R=49.5$ and $k\in \left[ 108.2,119.6\right].$
Thus, we obtain%
\begin{equation}
\left\{
	\begin{array}{ll}
		\left[ k/\left( 4\pi R\right) \right] ^{2} \leq 0.038 &\mbox{the case of one microsphere} \\
		\left[ 2k/\left( 4\pi R\right) \right] ^{2} \leq 0.150 &\mbox{the case of one microsphere} 
	\end{array}
\right.
\label{4.13}
\end{equation}

We observe that the number in the first line of (\ref{4.13}) almost
coincide with the bound in Figure \ref{fig 2a}. As to the numbers in the second
line of (\ref{4.13}), even though it overestimates somewhat the numerical bound in
Figure \ref{fig 2b}, still one can regard them as sufficiently small numbers.

\subsection{Reconstruction formulae}

\label{sec:4.3}

By \eqref{eqn usc}, for $k \in [\underline k, \overline k]$ and $\x \in \Pm$, we have
\begin{equation}
	|u|^2 = 1 + e^{\ik x_3} \overline {\usc} + e^{-\ik x_3} \usc + |\usc|^2.
	\label{4.1414}
\end{equation}
Using Theorem \ref{thm 4.1} and computational results of Figures \ref{fig 2a}--\ref{fig 2b},  
we drop the small term $\left\vert u_{\text{sc}}\right\vert ^{2}$ in (\ref{4.1414})
and obtain the following approximate formula: 
\begin{equation}
\left\vert u\left( \mathbf{x},k\right) \right\vert ^{2}=1+ 2\Re(e^{-\text{i}kx_{3}}%
u_{\text{sc}}),\quad\mathbf{x}%
\in P_{\text{meas}},k\in \left[ \underline{k},\overline{k}\right].  \label{4.15}
\end{equation}
Plugging \eqref{160} into \eqref{4.15}, we have
\begin{equation}
|u(\x, k)|^2 = -1 + 2A(\x)\Re \left(e^{\ik(\tau(\x) - x_3)}\right) = -1 + 2A(\x) \cos(k (\tau(\x) - x_3))
\label{4.1818}
\end{equation}
for all $k \in [\underline k, \overline k].$
Using Corollary \ref{corA}, we approximate the function $%
A\left( \mathbf{x}\right) $ as 
\begin{equation}
A(\mathbf{x})=|u(\mathbf{x},\overline{k})|=\sqrt{f(\mathbf{x},\overline{k})}%
,\quad \mathbf{x}\in P_{\text{meas}}.  \label{3.9}
\end{equation}%
Next, using using (\ref{4}), (\ref{4.1818}) and (\ref{3.9}), we
\begin{equation}
f\left( \mathbf{x}%
,k\right) =-1+2\sqrt{f(\mathbf{x},\overline{k})}\cos \left[ k\left( \tau
\left( \mathbf{x}\right) -x_{3}\right) \right] ,\quad \mathbf{x}\in P_{%
\text{meas}},k\in \left[ \underline{k},\overline{k}\right].\label{4.16}
\end{equation}
Hence, 
\begin{equation}
\tau \left( \mathbf{x}\right) =\frac{1}{k}\arccos \left( \frac{f\left( 
\mathbf{x},k\right) +1}{2\sqrt{f(\mathbf{x},\overline{k})}}\right) +R+\frac{%
2\pi m}{k}, \quad \mathbf{x}\in P_{%
\text{meas}},k\in \left[ \underline{k},\overline{k}\right]  \label{4.17}
\end{equation}
for some integer $m = m(\x, k).$ 
Plugging this function $\tau$ into \eqref{eqn asym} and noting that $e^{2\ii m \pi} = 1$, we obtain the following {\it inversion formula} 
\begin{equation}
u\left( \mathbf{x},k\right) =\sqrt{f(\mathbf{x},\overline{k})}\exp \left[ 
\text{i}\arccos \left( \frac{f\left( \mathbf{x},k\right) +1}{2\sqrt{f(%
\mathbf{x},\overline{k})}}\right) +\text{i}kR\right] ,\text{ }\mathbf{x}\in
P_{\text{meas}},k\in \left[ \underline{k},\overline{k}\right].  \label{4.18}
\end{equation}

\begin{remark}[Comment about the distance of measurement]\rm

\emph{\ It seems to be from both
(\ref{eqn LS}) and (\ref{4.5}) that, for a given value of the wave number }$%
k $\emph{, the larger the distance }$R$\emph{\ between the scatterers and
the measurement plane is, the better approximation of \ }$ u\left( 
\mathbf{x},k\right)$\emph{\ we would obtain by dropping the
term with }$\left\vert u_{\text{sc}}\right\vert ^{2}$ \emph{in (\ref{4.1414})}.
 \emph{However, this is not true} \emph{from both numerical
and Physics standpoints for exceedingly large values of }$R$\emph{.} \emph{%
Indeed, it follows from (\ref{4.3}) that, for a fixed value of }$k,$\emph{\
we have }$\left\vert u_{\text{sc}}\left( \mathbf{x},k\right) \right\vert
=O\left( R^{-1}\right) ,R\rightarrow \infty .$ \emph{Hence, by (\ref{4.1414})} 
$\left\vert u\left( \mathbf{x},k\right) \right\vert ^{2}=1+O\left(
R^{-1}\right) ,$ $R\rightarrow \infty .$\emph{\ Hence, for exceedingly large
values of }$R$ \emph{the influence of the terms with }$u_{\text{sc}}$\emph{\
in (\ref{4.1414}) is neglibly small, compared with $1.$} 
\emph{In other words, for these values of }$R$\emph{\ the total wave field
basically becomes the same as the incident plane wave is, i.e. without any
useful information about scatterers. This, therefore, makes it impossible to
reconstruct the complex valued wave field} $u\left( \mathbf{x},k\right) $ 
\emph{from the values of its intensity for those values of }$R$\emph{. }
\end{remark}

\begin{remark} \rm
\emph{Since formula (\ref{4.16}) is an approximate one and also since
the intensity data} $f\left( \mathbf{x},k\right) $ \emph{contains noise, 
the absolute value of the argument of the function }$\arccos $\emph{\ in (%
\ref{4.17}) and (\ref{4.18}) might exceed 1, at which we simply set in our computations $\tau
\left( \mathbf{x}\right) =R$.}
\end{remark}

\section{The phased coefficient inverse scattering problem}

\label{sec:5}

The reconstruction formula (\ref{4.18}) provides approximate values of the
wave field $u\left( \mathbf{x},k\right) $ at the measurement site $\mathbf{x}%
\in P_{\text{meas}}$ for the wave numbers $k\in \left[ \underline{k},%
\overline{k}\right] $ of our interval. This is done on the first stage of our
two-stage numerical procedure. We still need, however, to reconstruct the
unknown coefficient $n^{2}\left( \mathbf{x}\right) .$ And this is done on
the second stage of our procedure. Indeed, we have obtained now a phased
Coefficient Inverse Problem (phased CIP). It is well known that this
problem is not easy to solve since it is both nonlinear and ill-posed.
Still, we numerically solve this phased CIP on the second stage. We describe the
solution method of this problem in this section.

\begin{problem}[The phased coefficient inverse scattering problem]
Let $0<\underline{k}<\overline{k}<\infty $. Assume that the function $F(%
\mathbf{x},k),$ 
\begin{equation}
F(\mathbf{x},k)=u(\mathbf{x},k),\quad \mathbf{x}\in \Pm,k\in \lbrack 
\underline{k},\overline{k}]  \label{22}
\end{equation}%
is known. Determine the function $n(\mathbf{x})$ for $\mathbf{x}\in \Omega $%
. \label{Problem phased}
\end{problem}

Hence, by (\ref{21}), (\ref{3}) and (\ref{22}) 
\begin{equation}
F(\mathbf{x},k)=u(x_{1},x_{2},R,k), \quad \left\vert x_{1}\right\vert ,\left\vert
x_{2}\right\vert <b,k\in \left[ \underline{k},\overline{k}\right] ,
\label{23}
\end{equation}%
where the function $u(x_{1},x_{2},R,k)$ is defined in (\ref{4.18}).

Problem \ref{Problem phased} and its variations have been studied
extensively, see, e.g.~\cite%
{BeilinaKlibanovBook,ChowZou,KlibanovLiemLocHui:IPI2017} and references
cited therein. This problem is solved below by the globally convergent
numerical method, which was developed in \cite{KlibanovLiemLocHui:IPI2017}.
As it was mentioned in Section 1, this method was successfully tested on
microwave experimental data in \cite%
{AlekKlibanovLocLiemThanh:anm2017,LiemKlibanovLocAlekFiddyHui:jcp2017,nguyen2017:iip2017}%
.

Since the method of \cite{KlibanovLiemLocHui:IPI2017} was described in a
number of publications, we only briefly outline it below in this section for
the convenience of the reader. We refer to \cite{KlibanovLiemLocHui:IPI2017}
for details. The global convergence of this method is guaranteed by Theorem
6.1 of \cite{KlibanovLiemLocHui:IPI2017}.

We now comment on the issue of the uniqueness of this phased CIP. All
currently known uniqueness theorems for $n-$D ($n\geq 2)$ CIPs with single
measurement data are proven only for the case when the right hand side of
equation (\ref{eqn Helmholtz}) is not zero in $\overline{\Omega }.$ All such
theorems are proved using the idea of the paper \cite{BukhKlib}. In this
paper, the method of Carleman estimates was introduced in the field of
coefficient inverse problems. The idea of \cite{BukhKlib} became quite
popular since its inception, see, e.g. the most recent book \cite{BY},
sections 1.10 and 1.11 in the book \cite{BeilinaKlibanovBook}, the survey 
\cite{Kliba2013} and references cited therein. More recently the idea of 
\cite{BukhKlib} was extended to the construction of some globally convergent
numerical methods for CIPs, see e.g. \cite{Baud,KlibKam,Kliba2017b}.
However, the question about the uniqueness of the above phased CIP, so as
of some other similar ones, remains a well known open problem since the
right hand side of equation (\ref{eqn Helmholtz}) identically equals to
zero. Thus, we assume uniqueness for the sake of computations.

\subsection{Data propagation}

\label{sec:5.1}

The measurement plane $P$ is located far away from the targets. To avoid
working with a large computational domain, we \textquotedblleft move" the
data to a plane that is closer to the targets. More precisely, we move the
phased data to the plane containing the square $\Gamma \subset \partial
\Omega ,$ see (\ref{2}).\emph{\ }This procedure is called \textquotedblleft
data propagation." We briefly summarize it in Section \ref{sec:5.1}.

The data propagation aims approximately to \textquotedblleft move" the data
from the measurement plane to a plane that is close to the targets, named $%
P_p$ and $\Gamma \subset P_p.$ It was observed that this method also helps
to decrease the amount of noise in the data and focus the wave field on $P_p$%
. This method was extensively used for the preprocessing of experimental
data for the globally convergent algorithm, see \cite%
{AlekKlibanovLocLiemThanh:anm2017,
nguyen2017:iip2017,LiemKlibanovLocAlekFiddyHui:jcp2017}. We refer to~\cite%
{nguyen2017:iip2017} for a rigorous derivation of the formula~%
\eqref{propag_1} as well as for some more details about the data
propagation. The method is also known in optics as the \emph{spectrum
angular representation}, see \cite[Chapter 2]{NovotnyHecht:cup2012}.

Keeping in mind (\ref{23}), we define 
\begin{equation*}
\hat{u}(k_{1},k_{2},k)=\frac{1}{2\pi }\displaystyle \int\limits_{\mathbb{R}%
^{2}}F(x_{1},x_{2},k)\exp (-\mathrm{i}(k_{1}x_{1}+k_{2}x_{2}))dx_{1}dx_{2}.
\end{equation*}%
Here we extend by zero values of the function $F(x_{1},x_{2},k)$ for $%
\left\vert x_{1}\right\vert ,\left\vert x_{2}\right\vert >b.$ Finally the
propagated data is (see (\ref{1}), (\ref{2})) 
\begin{multline}
u(x_{1},x_{2},d_{2},k)=\frac{1}{2\pi }\displaystyle \int\limits_{%
\{k_{1}^{2}+k_{2}^{2}<k^{2}\}}\hat{u}(k_{1},k_{2},k)\exp \left( \mathrm{i}%
(k_{1}x_{1}+k_{2}x_{2}-k_{3}(d_{2}-R))\right) dk_{1}dk_{2},  \label{propag_1}
\end{multline}%
where 
\begin{equation}
k_{3}=\sqrt{k^{2}-k_{1}^{2}-k_{2}^{2}}, \quad \left\vert x_{1}\right\vert
,\left\vert x_{2}\right\vert <b\text{ and }k\in \left[ \underline{k},%
\overline{k}\right] .\text{ }  \label{24}
\end{equation}%
To speed up the process, in our computational implementation we compute the functions $\hat{u}%
(k_{1},k_{2},k)$ and $u(x_{1},x_{2},d_{2},k)$ using the Fast Fourier
transform. It follows from (\ref{2}) and (\ref{24}) that in (\ref{propag_1}), 
$\left( x_{1},x_{2},d_{2}\right) \in \Gamma .$ Hence, we set 
\begin{equation}
p\left( \mathbf{x},k\right) =u(x_{1},x_{2},d_{2},k)=u\mid _{\Gamma }, \quad k\in %
\left[ \underline{k},\overline{k}\right].  \label{25}
\end{equation}%
Since the boundary data $p\left( \mathbf{x},k\right) $ in (\ref{25}) are
given only on one side $\Gamma $ of the boundary of the domain $\Omega ,$ we
complement them on the rest of the boundary as 
\begin{equation}
\widetilde{p}\left( \mathbf{x},k\right) =\left\{ 
\begin{array}{ll}
p\left( \mathbf{x},k\right), &\mathbf{x\in }\Gamma , \\ 
e^{\text{i}kx_{3}}, &\mathbf{x\in }\partial \Omega \diagdown \Gamma .%
\end{array}%
\right.  \label{26}
\end{equation}%
Below we work only with the boundary data (\ref{26}). It was demonstrated
numerically in \cite{KlibanovLiemLocHui:IPI2017} that the reconstruction
result for the case of the complemented boundary data (\ref{26}) are close
to ones for the case when the computationally simulated data are assigned on
the entire boundary $\partial \Omega .$ The same conclusion was drawn in 
\cite{ChowZou} for the version of this method which is close to the one of 
\cite{BeilinaKlibanovBook}.\emph{\ }Also, (\ref{26}) was used in all above
cited works on microwave backscattering experimental data \cite%
{AlekKlibanovLocLiemThanh:anm2017,
nguyen2017:iip2017,LiemKlibanovLocAlekFiddyHui:jcp2017} and it did not negatively 
affect reconstruction results.

\subsection{An integro-differential equation}

\label{sec:5.2}

For applications in imaging of microscale objects, the wave numbers $%
\underline{k},\overline{k}$ are sufficiently large. Therefore, Theorem 3.1
implies that the function $u(\mathbf{x},k)$ is nonzero for $x\in \overline{%
\Omega },k\in \lbrack \underline{k},\overline{k}]$. Introduce the function $%
v(\mathbf{x},k)$ as 
\begin{equation}
v(\mathbf{x},k)=\log (u(\mathbf{x},k))\quad \mbox{for all }\mathbf{x}\in
\Omega ,k\in \lbrack \underline{k},\overline{k}].  \label{eqn log}
\end{equation}%
We refer to \cite[Section 4.1]{KlibanovLiemLocHui:IPI2017} for the
construction of the logarithm of our complex-valued function $u$. In any
case, since we use only derivatives of the function $v(\mathbf{x},k)$ rather
than this function itself, then the non-uniqueness of the logarithm of a
complex valued function is not of our concern. It is easy to see that 
\begin{equation}
\nabla v(\mathbf{x},k)=\frac{\nabla u(\mathbf{x},k)}{u(\mathbf{x},k)},\quad 
\mathbf{x}\in \Omega ,k\in \lbrack \underline{k},\overline{k}].
\label{varchange}
\end{equation}%
Using \eqref{eqn Helmholtz} and \eqref{varchange}, we obtain 
\begin{equation}
\Delta v(\mathbf{x},k)+(\nabla v(\mathbf{x},k))^{2}=-k^{2}n^{2}(\mathbf{x}).
\label{eqn v}
\end{equation}%
Differentiating the equation (\ref{eqn v}) with respect to $k$ and defining 
\begin{equation}
q(\mathbf{x},k)=\frac{\partial v(\mathbf{x},k)}{\partial k},  \label{q}
\end{equation}%
we obtain the following nonlinear integro-differential equation: 

\begin{equation}
\frac{k}{2}\Delta q(\mathbf{x},k)+k\nabla q(\mathbf{x},k)\Bigg(%
-\int\limits_{k}^{\overline{k}}\nabla q(\mathbf{x},s)ds+\nabla V(\mathbf{x})%
\Bigg)=  \label{eqn intdiff}
\end{equation}%
\begin{equation*}
-\int\limits_{k}^{\overline{k}}\Delta q(\mathbf{x},s)ds+\Delta V(\mathbf{x})+%
\Bigg(-\int\limits_{k}^{\overline{k}}\nabla q(\mathbf{x},s)ds+\nabla V(%
\mathbf{x})\Bigg)^{2},\quad \mathbf{x}\in \Omega ,k\in \lbrack \underline{k},%
\overline{k}],
\end{equation*}%
where the function 
\begin{equation}
V(\mathbf{x})=v(\mathbf{x},\overline{k}).  \label{eqntail}
\end{equation}%
is called the \textit{tail} function and it is unknown. The Dirichlet
boundary condition for the function $q(\mathbf{x},k)$ is%
\begin{equation}
q(\mathbf{x},k)=\frac{\partial _{k}\widetilde{p}\left( \mathbf{x},k\right) }{%
\widetilde{p}\left( \mathbf{x},k\right) },\quad \mathbf{x}\in \partial
\Omega ,k\in \lbrack \underline{k},\overline{k}],  \label{eqnBoundaryQ}
\end{equation}%
where the function $\widetilde{g}\left( \mathbf{x},k\right) $ is defined in (%
\ref{26}).

\subsection{The initial approximation for the tail function}

\label{sec:5.3}

The globally convergent numerical method is based on the iterative process
of solving equation (\ref{eqn intdiff}) with boundary condition (\ref%
{eqnBoundaryQ}) and simultaneously updating tail functions $V(\mathbf{x})$.
The whole procedure will be described in Algorithm \ref{algorithm phased}.
First, we explain in this section how to obtain the initial guess for the
tail function. We point out that when obtaining this guess, we do not use
any prior information about a small neighborhood of the exact solution of
our CIP.

Since the number $\overline{k}$ is assumed to be large, then, using (\ref{q}%
), (\ref{eqntail}), (\ref{eqn intdiff}) and the asymptotic behavior of the
function $v\left( \mathbf{x},k\right) $ at $k\rightarrow \infty ,$ which
follows from (\ref{eqn asym}), we obtain the following approximate equation
for the function $V(\mathbf{x})$ \cite[Section 4.1]%
{KlibanovLiemLocHui:IPI2017}: 
\begin{equation*}
\Delta V(\mathbf{x})=0,\quad \mathbf{x}\in \Omega .
\end{equation*}%
Note that to solve \eqref{eqn intdiff} and \eqref{eqnBoundaryQ}, we need to
know only the vector function $\nabla V(\mathbf{x})$. Thus, we can directly
compute the vector $\nabla V(\mathbf{x})$ via solution of the following
problem: 
\begin{equation}
\left\{ 
\begin{array}{rcll}
\Delta (\nabla V(\mathbf{x})) & = & 0 & \mbox{in }\Omega \\ 
\nabla V(\mathbf{x}) & = & R\left( \mathbf{x},\overline{k}\right) & 
\mbox{on
}\partial \Omega .%
\end{array}%
\right.  \label{eqn fa_problem}
\end{equation}%
Here we use (\ref{26}) to define the vector function $R\left( \mathbf{x},%
\overline{k}\right) =\left( R_{1},R_{2},R_{3}\right) \left( \mathbf{x},%
\overline{k}\right) $ as:%
\begin{equation}
R_{j}\left( \mathbf{x},\overline{k}\right) =\left\{ 
\begin{array}{ll}
\partial _{x_{j}}p\left( \mathbf{x},\overline{k}\right) /p\left( \mathbf{x},%
\overline{k}\right),& \mathbf{x}\in \Gamma , \\ 
0, & \mathbf{x}\in \partial \Omega \diagdown \Gamma ,%
\end{array}%
\right. j=1,2,  \label{27}
\end{equation}%
\begin{equation}
R_{3}\left( \mathbf{x},\overline{k}\right) =\left\{ 
\begin{array}{ll}
p_{1}\left( \mathbf{x},\overline{k}\right) /p\left( \mathbf{x},\overline{k}%
\right), &\mathbf{x}\in \Gamma , \\ 
i\overline{k}, &\mathbf{x}\in \partial \Omega \diagdown \Gamma .%
\end{array}%
\right.  \label{28}
\end{equation}%
Here $p_{1}\left( \mathbf{x},\overline{k}\right) =\partial
_{x_{3}}u(x_{1},x_{2},d_{2},\overline{k})=\partial _{x_{3}}u\left( \mathbf{x}%
,\overline{k}\right) \mid _{\Gamma }.$ To find the function $p_{1}\left( 
\mathbf{x},\overline{k}\right) ,$ we propagate the data as in Section 6.1 to
two planes: $P=\left\{ x_{3}=d_{2}\right\} $ as in (\ref{propag_1}) and $%
P^{\varepsilon }=\left\{ x_{3}=d_{2}+\varepsilon \right\} $ for a
sufficiently small number $\varepsilon >0.$ Next, we use 
\begin{equation}
p_{1}\left( \mathbf{x},\overline{k}\right) \approx \frac{%
u(x_{1},x_{2},d_{2}+\varepsilon ,\overline{k})-u(x_{1},x_{2},d_{2},\overline{%
k})}{\varepsilon }.  \label{29}
\end{equation}

The solution of problem \eqref{eqn fa_problem} is considered as an initial
approximation of $\nabla V$, which is denoted by $\nabla V_{0}(\mathbf{x})$.
It was also derived in \cite{KlibanovLiemLocHui:IPI2017} from the above
mentioned asymptotic behavior of the function $v\left( \mathbf{x},k\right) $
that the approximate value of the vector function $\nabla q\left( \mathbf{x},%
\overline{k}\right) $ can be found as%
\begin{equation}
\nabla q\left( \mathbf{x},\overline{k}\right) =\frac{\nabla V_{0}(\mathbf{x})%
}{\overline{k}}.  \label{30}
\end{equation}

\begin{remark}\rm
\begin{enumerate}
\item \emph{It was shown numerically in \cite{KlibanovLiemLocHui:IPI2017}
that formulae (\ref{27}) and (\ref{28}) provide almost the same results as
ones in the case when correct boundary conditions are prescribed on the
entire boundary} $\partial \Omega .$

\item \emph{It follows from (\ref{eqnBoundaryQ}) and (\ref{29}) that we
should differentiate the data, which usually contain noise. In all such
cases we use finite differences for the differentiation. However, neither in
the current work, nor in all above cited works with noisy experimental data
we have not observed any instabilities, probably because the step sizes of
our finite differences were not too small.}
\end{enumerate}
\end{remark}
\subsection{The globally convergent algorithm for Problem \protect\ref%
{Problem phased}}

\label{sec:5.4}

As soon as the first approximation $\nabla V_{0}(\mathbf{x})$ for the
gradient of the function $V(\mathbf{x})$ is found, we can proceed with our
iterative algorithm of solving the problem (\ref{eqn intdiff}) and (\ref%
{eqnBoundaryQ}). So, we present in this section the globally convergent
algorithm which we use in this paper. Let 
\begin{equation*}
k_{N}=\underline{k}<k_{N-1}<k_{N-2}<\dots <k_{1}<k_{0}=\overline{k}
\end{equation*}%
be a uniform partition of $[\underline{k},\overline{k}]$. Let $%
h=k_{i}-k_{i+1}$, $i=0,...,N-1$ denotes the grid step size of this
partition. Define the function $q_{m}(\mathbf{x})$ as 
\begin{equation*}
q_{m}(\mathbf{x})=q(\mathbf{x},k_{m}),\quad u_{m}(\mathbf{x})=u(\mathbf{x}%
,k_{m}),\quad \mathbf{x}\in \Omega .
\end{equation*}%
Let the function $\nabla V_{0}$ be calculated as in Section \ref{sec:5.3}
and the vector function $\nabla q_{0}(\mathbf{x})=\nabla q(\mathbf{x},%
\overline{k})$ is approximately found as in (\ref{30}). Using the
mathematical induction, assume that vector functions $\nabla V_{n-1}$ and $%
\nabla q_{n-1}$ are known for some $m\in \{1,\dots ,N\}$. Recall that $%
\Delta \varphi =\func{div}\left( \nabla \varphi \right) $ for any
appropriate function $\varphi \left( \mathbf{x}\right) $. Thus,
integro-differential equation (\ref{eqn intdiff}) is approximated by 
\begin{equation*}
\frac{k_{m}}{2}\Delta q_{m}(\mathbf{x})+k_{m}\nabla q_{m}(\mathbf{x})\nabla
Q_{m-1}(\mathbf{x})+k_{m}\nabla q_{m-1}(\mathbf{x})\nabla V_{m-1}(\mathbf{x}%
)=
\end{equation*}%
\begin{equation}
-\Delta Q_{m-1}(\mathbf{x})+\Delta V_{m-1}(\mathbf{x})+\big(-\nabla Q_{m-1}(%
\mathbf{x})+\nabla V_{m-1}(\mathbf{x})\big)^{2},\quad \mathbf{x}\in \Omega ,
\label{bvp}
\end{equation}%
where 
\begin{equation*}
\nabla Q_{m-1}(\mathbf{x})=h\sum\limits_{j=0}^{m-1}\nabla q_{j}(\mathbf{x}),
\end{equation*}%
is an approximation of the integral%
\begin{equation*}
\int_{k}^{\overline{k}}\nabla q(\mathbf{x},s)ds.
\end{equation*}%
By (\ref{eqnBoundaryQ}) the boundary condition for the function $q_{n}(%
\mathbf{x})$ can be calculated as 
\begin{equation}
q_{n}(\mathbf{x})=\frac{\partial _{k}\widetilde{p}\left( \mathbf{x}%
,k_{n}\right) }{\widetilde{p}\left( \mathbf{x},k_{n}\right) },\quad \mathbf{x%
}\in \partial \Omega .  \label{bvp boundary}
\end{equation}

Our globally convergent algorithm of the solution of\ Problem \ref{Problem
phased} is summarized in Algorithm \ref{algorithm phased}.

\begin{algorithm}
	\caption{Globally convergent algorithm for Problem \ref{Problem phased}.}
	\label{algorithm phased}
	\begin{algorithmic}[1]
		\State{\label{stepPropagation} Propagate the data from $\Pm$ to $\Pp$ as in Section \ref{sec:5.1}}.
		\State{Calculate the 
			initial approximation $\nabla V_0(\x)$ of $\nabla V(\x)$ as in Section \ref{sec:5.3}.  Set $\nabla q_0(\x)$ as in \eqref{30}.}
		\For{$n=1,2,\dots ,N$}
			\State{Assume, by induction, that $\nabla V_{m}(\x)$ and $\nabla q_{m}(\x)$, $m = 0, \dots, n - 1$, are known. Set \[ Q_{n - 1}(\x) = h\sum_{m = 0}^{n - 1}q_{m}(\x).\] Calculate $\nabla V_n(\x)$ and $\nabla q_n(\x)$ as follows.}
			\State{Set $\nabla q_{n,0}(\x) = \nabla q_{n-1}(\x)$ and $\nabla V_{n,0}(\x) = \nabla V_{n-1}(\x)$.}
			\For{$i=1,2,\dots ,I_N$} 
				\State{Find $q_{n,i}(\x)$ by solving the boundary value problem 
					\eqref{bvp}-\eqref{bvp boundary}.}
				\State{Update $\nabla v_{n,i}(\x) =\ds -\left(h\nabla q_{n,i}(\x) +Q_{n - 1}(\x) \right) + \nabla V_{n,i-1}(\x)$. }
				\State{\label{step cni} Update $c_{n,i}(\x)$ via \eqref{eqn v}.}
				\State{Find $u_{n,i}(\mathbf{x},\overline k)$ by solving integral equation 
					\eqref{eqn LS} with $\beta(\mathbf{x}) = c_{n,i}(\mathbf{x}) -1$.}
				\State{Update $\nabla V_{n,i}(\mathbf{x}) = \nabla u_{n,i}
					(\mathbf{x},\overline k)/u_{n,i}(\mathbf{x},\overline k)$.}
			\EndFor
			\State{Update $q_n(\x) = q_{n,I_N}(\x), c_n = c_{n,I_N}(\x)$ and $\nabla V_n(\x) = \nabla V_{n,I_N}(\x)$.}
		\EndFor
	\State{The function $c(\x)$ is set to be the function $c_{n_*}(\x)$, in which
	\[n_* = \mbox{argmin}\left\{\frac{\|c_{n-1} - c_n\|}{\|c_n\|}: n = 3, \dots N, \right\}.\]
	}
	\end{algorithmic}
\end{algorithm}

The Dirichlet boundary value problem (\ref{bvp}), (\ref{bvp boundary}) is
solved via the finite element method using the software FreeFem++ (see~\cite%
{Hecht:jmm2012}). The numerical solution of the Lippmann-Schwinger integral
equation in Step 11 is performed using the numerical method developed in 
\cite{LechleiterNguyen:acm2014, Nguyen:anm2015}. The number $I_{N}$ of inner
iterations is chosen computationally. Usually $I_{N}=3$, see details in \cite%
{LiemKlibanovLocAlekFiddyHui:jcp2017}.

\section{Data collection}

\label{sec:6}

\begin{figure}[tbp]
\begin{center}
	\subfloat[\label{fig 3a} \emph{The photograph of our experimental device.}]{
		\includegraphics[width=0.35\linewidth]{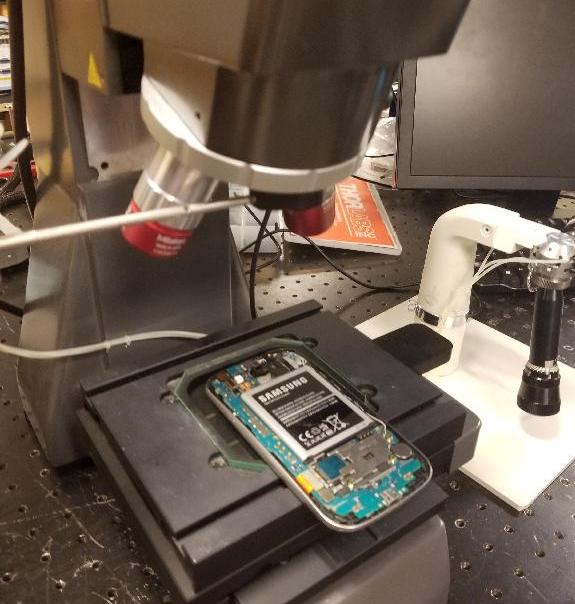}}
	\subfloat[\label{fig 3b} \emph{A schematic diagram of data collection.}]{
		\includegraphics[width=0.5\linewidth]{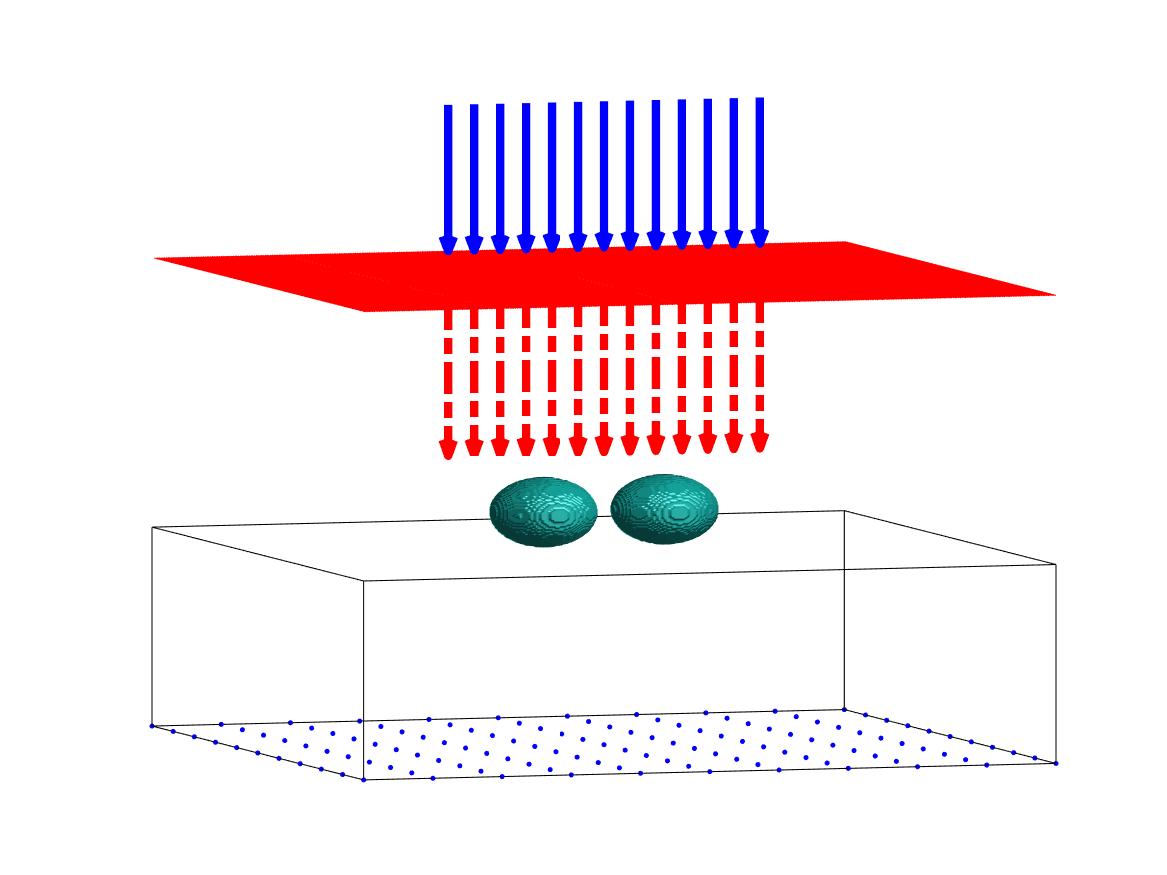}
		\put(-200,20){\tiny $P_{meas}$}	
		\put(-110,11){\tiny Detectors}
		\put(-145,40){\tiny Protective glass layer}
		\put(-80, 75){\tiny Microspheres}
		\put(-70, 95){\tiny Filtered light}
		\put(-120, 115){\tiny Filter}
		\put(-70, 140){\tiny White light}
		}
\end{center}
	\caption{\emph{The data collection setup}}
	\label{fig:datacol}
\end{figure}
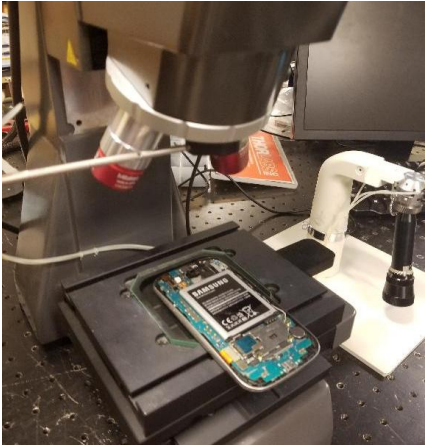
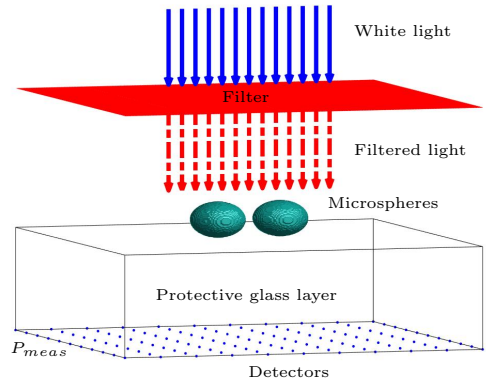

As it was stated in Introduction, following the idea of a series of works of
experimentalists, such as, e.g. \cite{Tseng}, we use the detector array of
the camera of Samsung Galaxy S3 cell phone. Figure \ref{fig:datacol} displays the photograph of our
experimental device. We use white light. By using narrow-pass spectral filters, 
we filtered relatively narrow bands with 
the width around 10$\mu m$ centered at the different wavelengths ($\lambda$) which 
were distributed throughout the visible spectrum: 
\begin{equation}
\lambda = 0420, 0.473, 0.525, 0.580, 0.620, 
0.671 \mu m.
\label{6.1}
\end{equation}
We have made variables dimensionless via the change of variables in
Helmholtz equation (\ref{eqn Helmholtz}) 
\begin{equation}
\mathbf{x}^{\prime }=\frac{\mathbf{x}}{10\mu m}.  \label{6.10}
\end{equation}%
Then the dimensionless wavelength $\lambda ^{\prime }=\lambda /(10\mu m)$ and the
wave number $k=2\pi /\lambda ^{\prime }.$We keep the same notations for new
variables for brevity. Hence, by (\ref{6.1}) we use the following values of $%
k$:%
\begin{equation}
k=93.6,101.3,108.3,119.7,132.8,149.5.  \label{6.2}
\end{equation}%
However, as stated in Section \ref{sec:4}, we have observed that the data
is too noisy for all values of $k$ in (\ref{6.2}). Only data for wave numbers $108.3$ and $%
119.7$ have acceptable noise level. Hence, we set as in (\ref{4.1}) $\underline{k}=108.3,\overline{k}%
=119.7$. Next, to obtain the data for $\left\vert u\left( \mathbf{x}%
,k\right) \right\vert ^{2}\mid _{\mathbf{x}\in P_{\text{meas}}}$ for all
other values of $k\in \left[ \underline{k},\overline{k}\right] ,$ we have
linearly interpolated the measured data between points $k=\underline{k}$ and 
$k=\overline{k}.$

An array of camera's photosensitive pixels was covered by the manufacturer 
with a protective glass layer, as illustrated in Figure \ref{fig 3b}.
In this glass, the refractive index $n_{0}=1.5$.
Thus, we take this value as the refractive index of the background. We must
scale the model, so that the value of the refractive index in the background
would become $n_{\text{bkgr}}=1$. To do this, we again change variables as 
\begin{equation}
\mathbf{y}=n_{0}\mathbf{x},\text{ }\mathbf{x}\in \mathbb{R}^{3}.  \label{6.3}
\end{equation}%
The Helmholtz equation in (\ref{eqn Helmholtz}) becomes 
\begin{equation*}
\Delta _{\mathbf{y}}u(\mathbf{y},k)+k^{2}\left( \frac{n(\mathbf{y})}{n_{0}}%
\right) ^{2}u(\mathbf{y},k)=0.
\end{equation*}%
When we solve Problem \ref{Problem phased} by the above globally convergent
method, we find first the relative contrast $n_{\text{rel}},$ 
\begin{equation*}
n_{\text{rel}}=\frac{n}{n_{0}}=\frac{n}{1.5}.\text{ }
\end{equation*}%
Next, we find the refractive index as $n=1.5n_{\text{rel}}.$

Below all sizes are those which were made dimensionless first by (\ref{6.10}%
) and scaled ones then by (\ref{6.3}). The value of the dielectric constant in
each microsphere of our experiments was 
\begin{equation}
n\left( \text{microsphere}\right) =2.15,\text{ }n_{\text{rel}}\left( \text{%
microsphere}\right) =1.43.  \label{6.5}
\end{equation}

In our experiments, we have collected the above mentioned data for the case
when scatterers were microspheres of the radius $0.45$. The center of each
microsphere was located on the plane $\left\{ x_{3}=0\right\} .$ The
measurement plane was located at $R=49.5$. We have measured the data at that
plane on a square $P_{\text{meas}},$ 
\begin{equation}
P_{\text{meas}}=\left\{ \left\vert x_{1}\right\vert ,\left\vert
x_{2}\right\vert <3.75,x_{3}=R=49.5\right\} ,  \label{6.6}
\end{equation}%
\begin{equation*}
\Omega =\left\{ \left\vert x_{1}\right\vert ,\left\vert x_{2}\right\vert
<3.75,x_{3}\in \left( 6.8,0.7\right) \right\} ,
\end{equation*}%
\begin{equation}
\Gamma =\left\{ \left\vert x_{1}\right\vert ,\left\vert x_{2}\right\vert
<3.75,x_{3}=0.7\right\} .  \label{6.8}
\end{equation}

\section{Numerical Results{\ }}

\label{sec:7}

In this section, we present our numerical results for solving Problem \ref%
{Problem phaseless} for both computationally simulated and experimental
data. First, we have reconstructed the wave field at $P_{\text{meas}}$ using
the reconstruction formula (\ref{4.18}). 
We reduce  Problem \ref{Problem phaseless} to Problem \ref{Problem phased}. Next, we
have solved the latter problem by the method outlined in Section \ref{sec:5}
and have reconstructed the unknown refractive index $n\left( \mathbf{x}%
\right) $ this way, i.e. we have imaged our microspheres.

We had two sets of experimental data. The first one was for the case of a
single microsphere and the second one was for the case of two microspheres.
For each case we have solved both phaseless CIPs: one with experimental data and the
second one with computationally simulated data. It is important that
computationally simulated data were generated for exactly the same values of
parameters as those in corresponding microspheres in experiments: we have
used exactly the same radii of spheres $0.45$, locations of their centers
and values of the relative index $n=2.15$ in (\ref{6.5}). In all tests, the
interval of wave numbers was the sama as in (\ref{4.1}).

We have performed such tests for computationally simulated data in order to
verify our method via the comparison of results with those for experimental
data.

We now describe how we modeled our microspheres for computationally
simulated data. In order to improve the accuracy of the solution of the
forward problem, we have decided to smooth out inclusion/background
interface rather than having the function $c\left( \mathbf{x}\right) $ which
would have a discontinuity on this interface.

Define the smooth function $\psi \in C^{\infty }(\mathbb{R}^{3})$ as 
\begin{equation*}
\psi (\mathbf{x})=\left\{ 
\begin{array}{ll}
\exp \left( -\frac{|\mathbf{x}|^{2}}{1-|\mathbf{x}|^{2}}\right) & |\mathbf{x}%
|<1 \\ 
0 & |\mathbf{x}|\geq 1.%
\end{array}%
\right.
\end{equation*}%
Let $\mathbf{x}_{0}$ be the center and $r=0.45$ be the radius of the
microsphere used in our simulations. Then, taking into account that by (\ref%
{6.5}) $n_{\text{rel}}^{2}\left( \text{microsphere}\right) =2.04$, we use
the following function $n^{2}(\mathbf{x})$ in our computational simulations%
\begin{equation*}
n^{2}(\mathbf{x})=1+1.04\psi \left( \frac{\mathbf{x}-\mathbf{x}_{0}}{r}%
\right) .
\end{equation*}%
In the case of two microspheres with their centers at $\mathbf{x}%
_{0}^{\left( 1\right) }$ and $\mathbf{x}_{0}^{\left( 2\right) }$ and with
the distance between their centers exceeding 1, we have%
\begin{equation*}
n^{2}(\mathbf{x})=1+1.04\psi \left( \frac{\mathbf{x}-\mathbf{x}_{0}^{\left(
1\right) }}{r}\right) +1.04\psi \left( \frac{\mathbf{x}-\mathbf{x}%
_{0}^{\left( 2\right) }}{r}\right) .
\end{equation*}%
Let $n^{\text{comp}}\left( \mathbf{x}\right) ,\mathbf{x}\in \overline{\Omega 
}$ be the computed refractive index $n\left( \mathbf{x}\right) .$ As the
computed value $n_{\text{comp}}$ of the refractive index $n$ in inclusions,
we take:%
\begin{equation}
n_{\text{comp}}\left( \text{microsphere}\right) =\max_{\overline{\Omega }}n^{%
\text{comp}}\left( \mathbf{x}\right) .  \label{6.9}
\end{equation}

\subsection{Test 1: One microsphere}

\label{sec:7.1}

\begin{figure}[tbp]
	\subfloat[\it Experimental data.]{\includegraphics[width=0.5\linewidth]{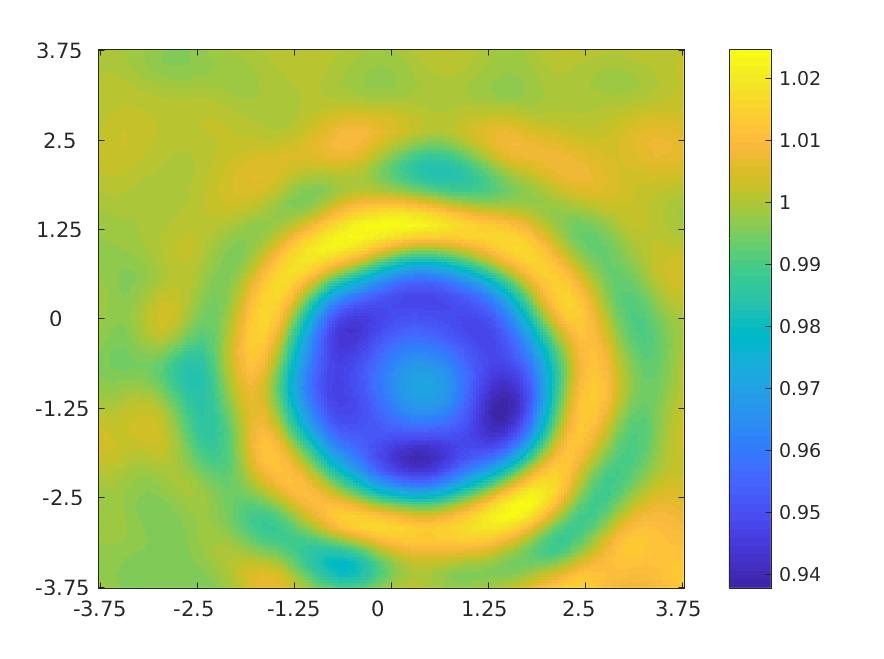} }
	\subfloat[\it Computationally simulated data.]{\includegraphics[width=0.5\linewidth]{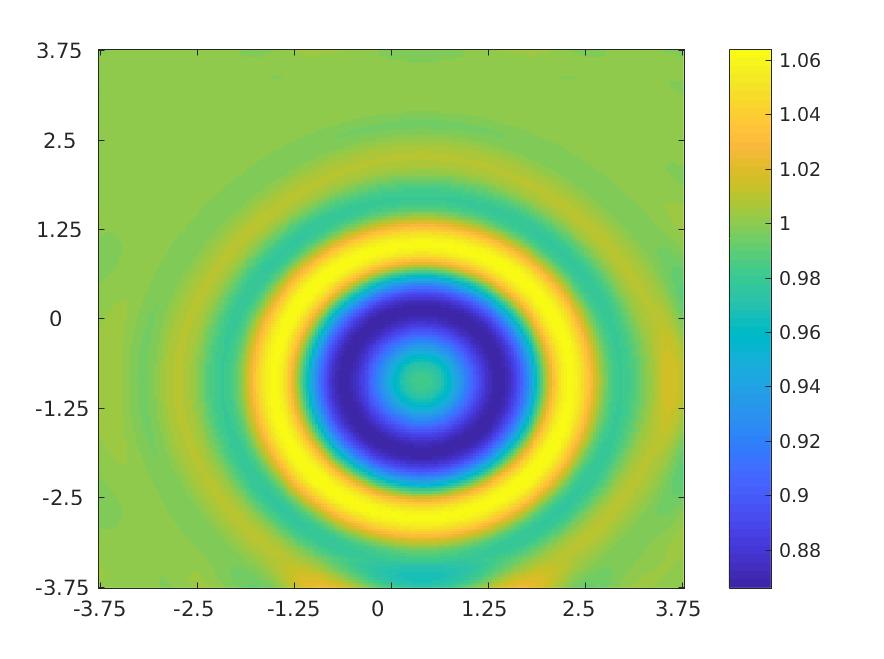} }
	\caption{\label{fig 4} \it Test 1. The function $|u(\x, k)|$, $\x \in \Pm$ and $k = 119.7$. The similarity of experimental data and the computationally simulated data is evident.}
\end{figure}

Figure \ref{fig 4} displays the modulus $\left\vert u\left( \mathbf{x},\underline{k}%
\right) \right\vert ,\underline{k}=108.3$ of both experimental (a) and
computationally simulated (b) data at the measurement square $P_{\text{meas}%
},$ see (\ref{6.6}).\ One can see good agreement between the experimental (left) 
and computed (right) data. The difference can be ascribed to the noise in the experimental
system.

\begin{figure}[tbp]
	\subfloat[\it The true and reconstructed real parts of the function $u(\x, k)$]{\includegraphics[width=0.45\linewidth]{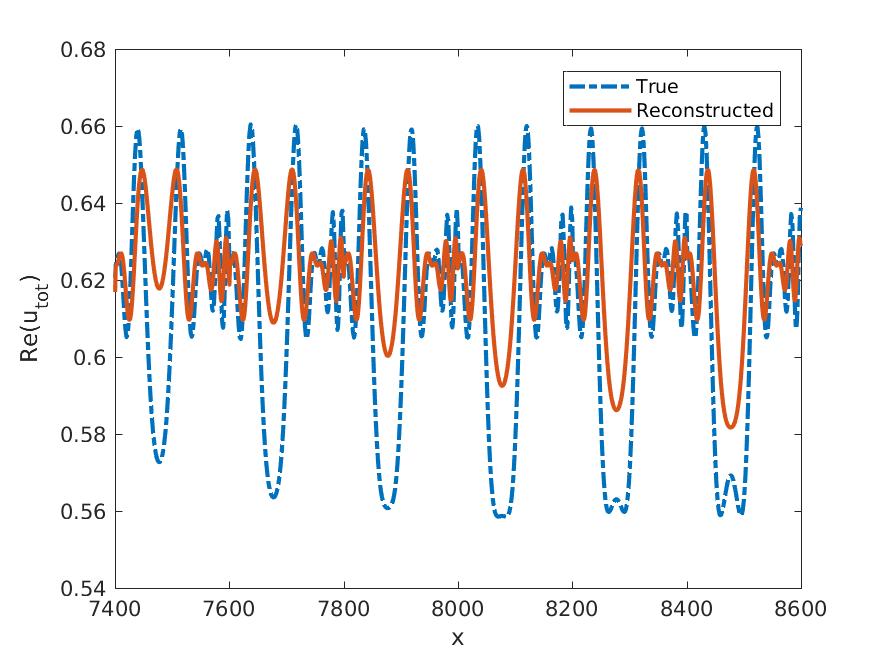}}\hfill
	\subfloat[\it The true and reconstructed imaginary parts of the function $u(\x, k)$] {\includegraphics[width=0.45\linewidth]{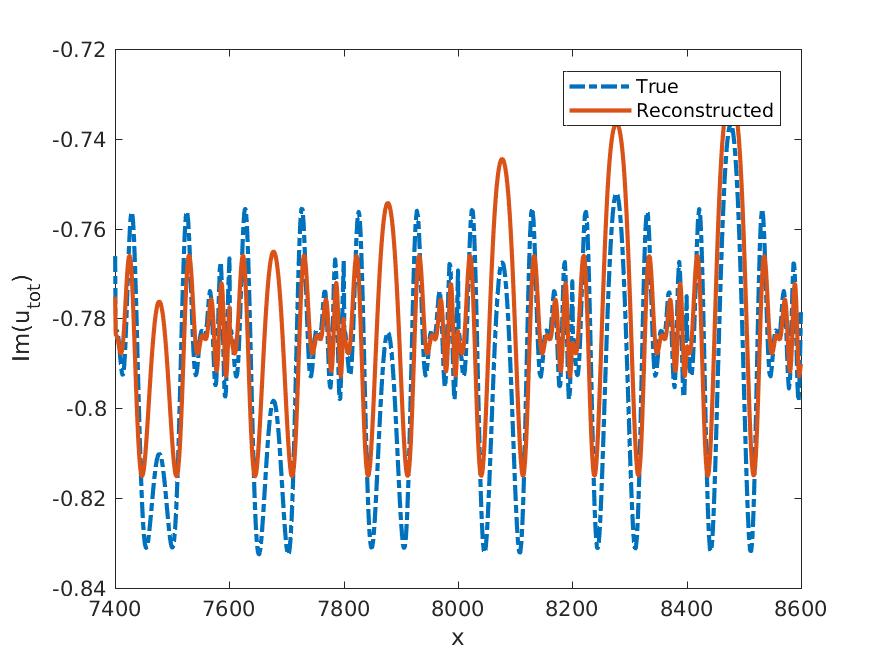}}
	\caption{\label{fig 5} \it Test 1. The reconstruction of the complex valued function $u(\x, k)$ with $k = 119.7$ using the inversion formula \eqref{4.18}. 
	One can see that imaginary parts coincide quite well, whereas real parts coincide sort 
	of satisfactory, given a highly oscillatory behavior of these curves. We conclude that 
	our reconstruction formula \eqref{4.18} is quite accurate. }
\end{figure}

Figure \ref{fig 5} represents the true and reconstructed real and imaginary parts of the
function $u\left( \mathbf{x},\underline{k}\right) $ for $\mathbf{x}\in P_{%
\text{meas}}$ for computationally simulated data.\ The true function $%
u\left( \mathbf{x},\underline{k}\right) $ was computed via the numerical
solution of equation (\ref{eqn LS}). The reconstructed function $u\left( 
\mathbf{x},\underline{k}\right) $ was computed via the reconstruction
formula (\ref{4.18}). Figures 6 were obtained as follows: we have arranged a
uniform $100\times 100$ grid of points $\{\mathbf{x}_{l}\}_{l=1}^{10,000}$
in $P_{\mathrm{meas}}$. So, Figures 6 show true and reconstructed values of $%
\func{Re}u(\mathbf{x}_{s},\underline{k})$ and $\func{Im}u(\mathbf{x}_{s},%
\underline{k})$, where $s\in \{7400,\dots ,8600\}$. One can see that
imaginary parts coincide quite well, whereas real parts coincide sort of
satisfactory, given a highly oscillatory behavior of these curves. In
addition, using Remark 4.1, we conclude that our reconstruction formulae of
Theorem 4.2 are rather accurate ones, especially given a difficult nature of
our phaseless CIP.

\begin{figure}[tbp]
\subfloat[\label{fig 6a}\it The propagation of experimental data]{\includegraphics[width=0.5\linewidth]{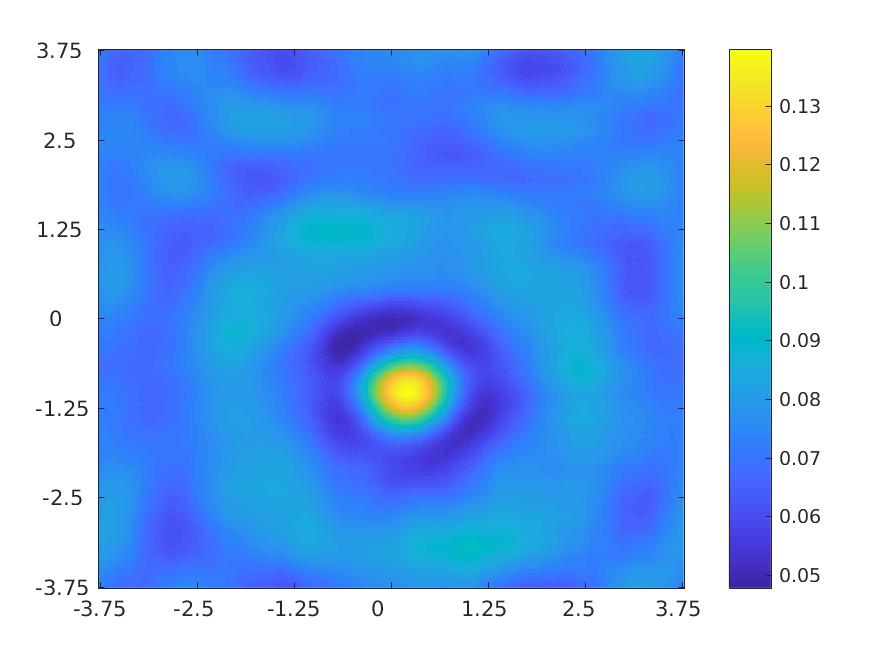}}
\subfloat[\label{fig 6b}\it The propagation of simulated data]{\includegraphics[width=0.5\linewidth]{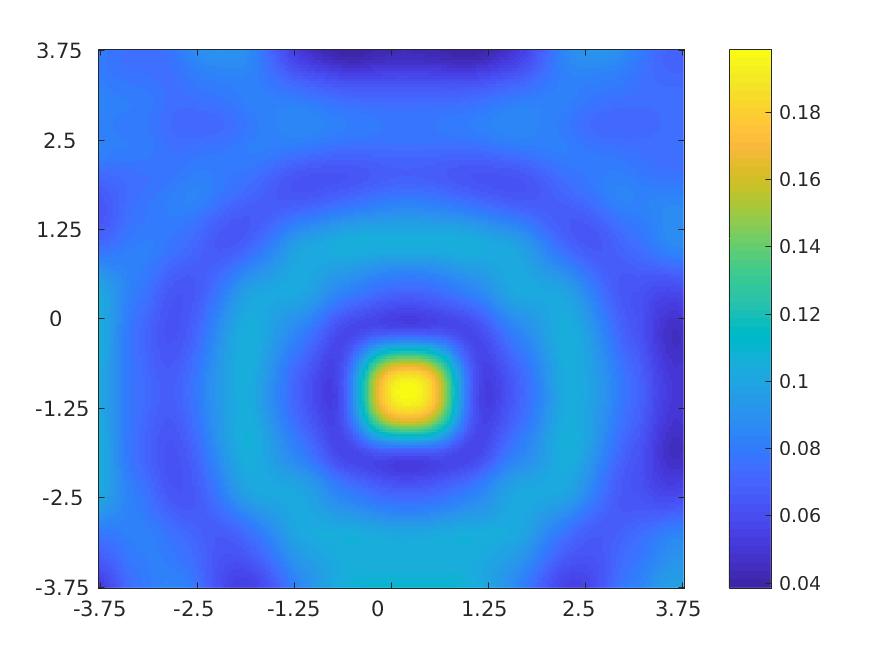}}
\caption{\label{fig 6}\it Test 1. Modulus of the propagated data $|u (\x, k)|$ for $\x \in \Gamma$.  A good similarity between (a) and (b) is observed. One can also see that the data propagation helps to figure out $(x_1, x_2)$ coordinates of inclusions.}
\end{figure}

Figure \ref{fig 6} displays the modulus of $\left\vert u\left( \mathbf{x},\underline{k}%
\right) \right\vert $ for $\mathbf{x}\in \Gamma ,$ where $\Gamma \subset
\partial \Omega $ was defined in (\ref{6.8}). In other words, Figures \ref{fig 6} show
the modulus of propagated data, both for the experimental and
computationally simulated cases. One can see that Figures \ref{fig 6a} and \ref{fig 6b} look
quite similar. In other words, propagated data for these two cases are quite
close to each other. One can also see that the data propagation helps to
figure out $\left( x_{1},x_{2}\right) $ coordinates of inclusions.

\begin{table}[tbp]
\caption{\label{table 1} Test 1. True and computed values of the refractive indices
in reconstructed microspheres, see (\protect\ref{6.9}).}
\begin{center}
\begin{tabular}{|l|l|l|l|}
\hline
Data & True $n\left( \text{microsphere}\right) $ & $n_{\text{comp}}\left( 
\text{microsphere}\right) $ & Relative error \\ \hline
Experimental & 2.15 & 2.043 & 5.2\% \\ \hline
Simulated \  & 2.15 & 2.04 & 5.4\% \\ \hline
\end{tabular}%
\end{center}
\end{table}

It follows from Table \ref{table 1} that values of the computed refractive indices of
reconstructed microspheres for experimental and simulated data are almost
the same. The computational error is quite small in both cases.

\begin{figure}[tbp]
	\subfloat[\it True image]{\includegraphics[width=0.3\linewidth]{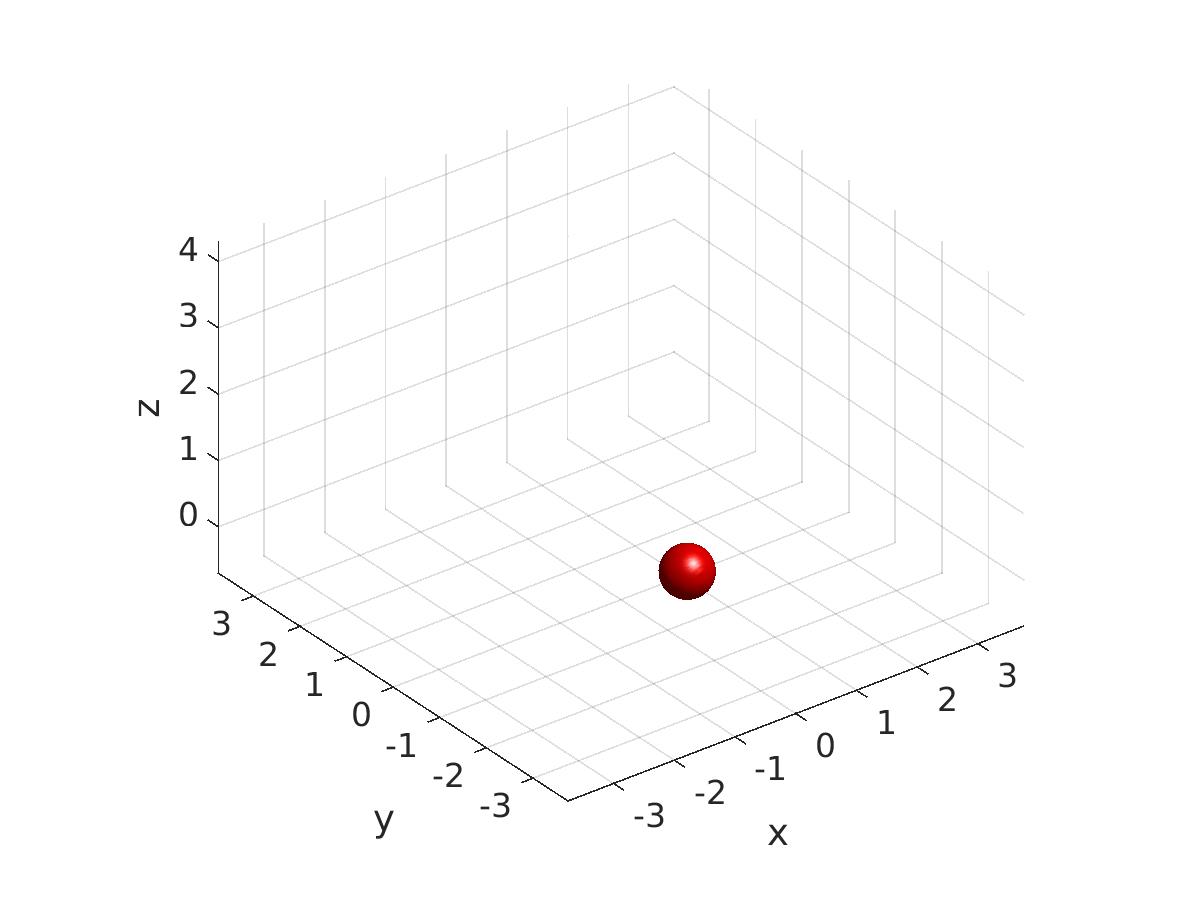}}\hfill
	\subfloat[\it Image computed from simulated data.]{\includegraphics[width=0.3\linewidth]{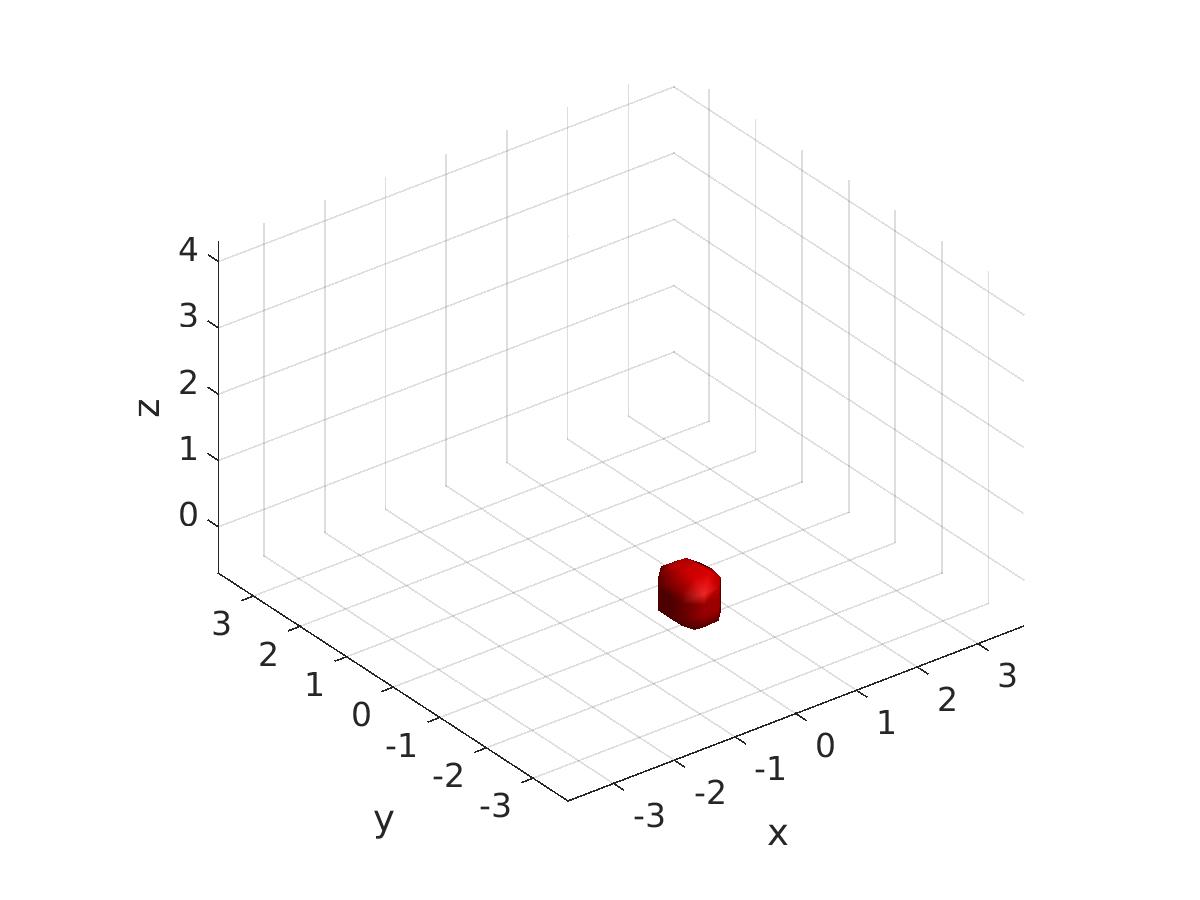}}\hfill
	\subfloat[\it Image computed from experimental data.]{\includegraphics[width=0.3\linewidth]{Figure8b}}
	\caption{\label{fig 7} \it Test 1. True and computed images of the microsphere involved in the experiment. One can observe that the location of the unknown microsphere is imaged with a very good accuracy in both cases. A similarity between (b) and (c) can also be observed.}
\end{figure}

Figure \ref{fig 7} displays resulting images. One can see that images for both
computationally simulated and experimental data are very similar and also
similar with the true image. Thus, locations of inclusions are also
correctly reconstructed.

\subsection{Test 2: Two microspheres}

\label{sec:7.2}

Since results for this case are very similar with ones for the case of one
microsphere, we shorten in this section, compared with the previous one.
Since the comparison of functions $\func{Re}u(\mathbf{x},\underline{k})$ and 
$\func{Im}u(\mathbf{x},\underline{k})$ for $\mathbf{x}\in P_{\text{meas}}$
for computationally simulated and experimental data is very similar for this
case with Test 1, we do not show here an analog of Figure \ref{fig 5}: for brevity.
The same for an analog of Figure \ref{fig 6}. %

\begin{figure}[tbp]
	\subfloat[\it Experimental data]{\includegraphics[width=0.45\linewidth]{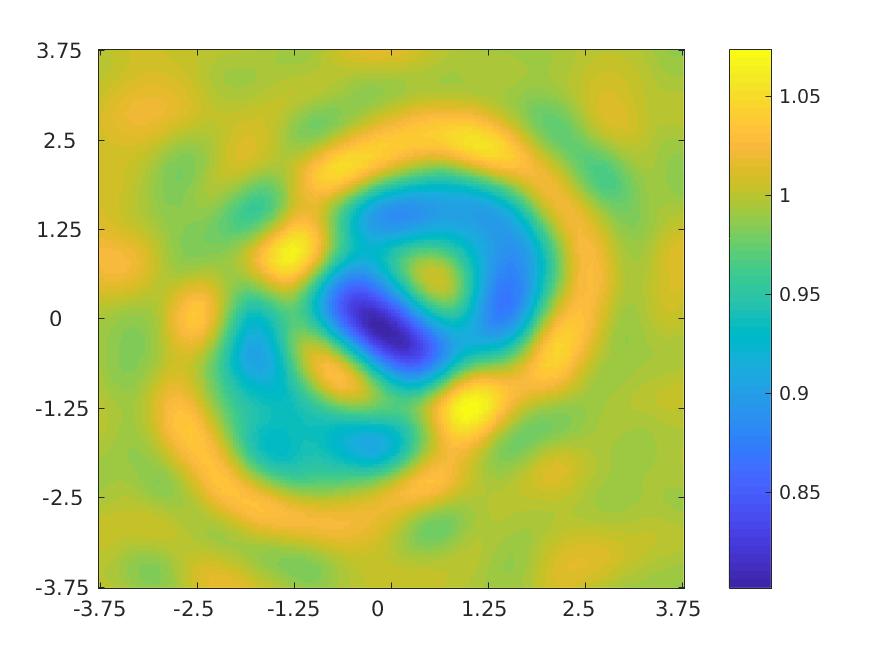}} \hfill
	\subfloat[Simulated data]{\includegraphics[width=0.45\linewidth]{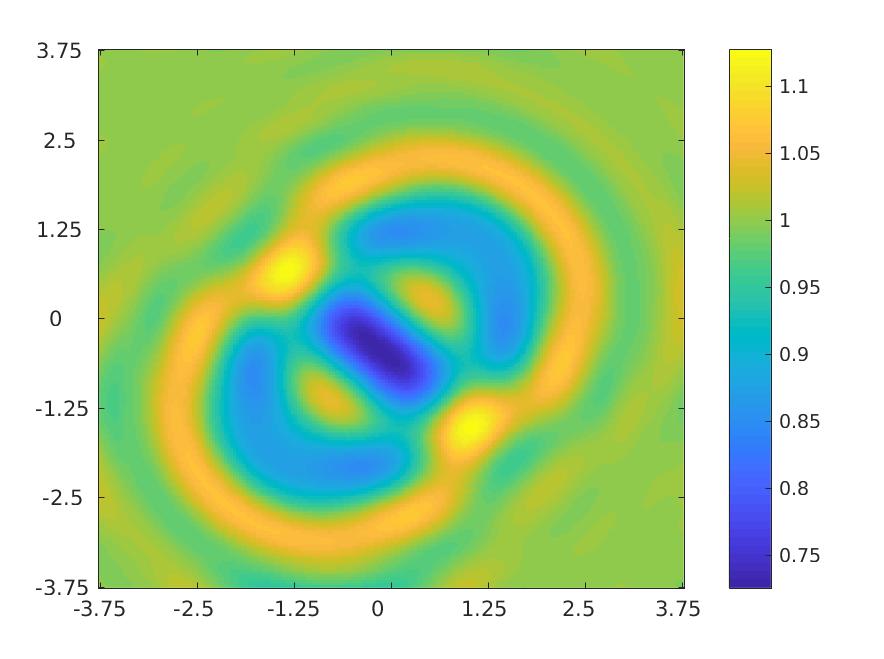}}
	\caption{\label{fig 8}\it Test 2. The data $|u(\x, k)|$, $\x \in \Pm$ and $k = 119.7$.}
\end{figure}

\begin{table}[tbp]
\caption{\label{table 2} Test 2. Correct and computed values of the refractive indices in
two reconstructed microspheres, see (\protect\ref{6.9}).}
\begin{center}
\begin{tabular}{|l|l|l|l|}
\hline
Data & Correct $n\left( \text{microsphere}\right) $ & $n_{\text{comp}}\left( 
\text{microsphere}\right) $ & Relative error \\ \hline
Experimental & 2.15 & 2.04 (in both microspheres) & 5.4\% \\ \hline
Simulated \  & 2.15 & 2.04 (in both microspheres) & 5.4\% \\ \hline
\end{tabular}%
\end{center}
\end{table}

\begin{figure}[tbp]
	\subfloat[\it True image]{\includegraphics[width=0.3\linewidth]{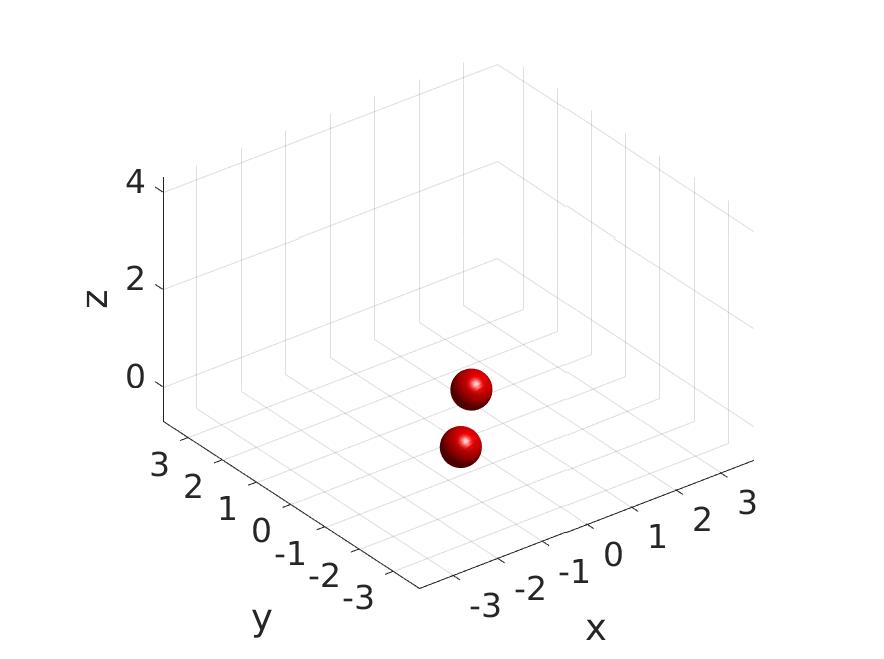}}\hfill
	\subfloat[\it Image computed from simulated data]{\includegraphics[width=0.3\linewidth]{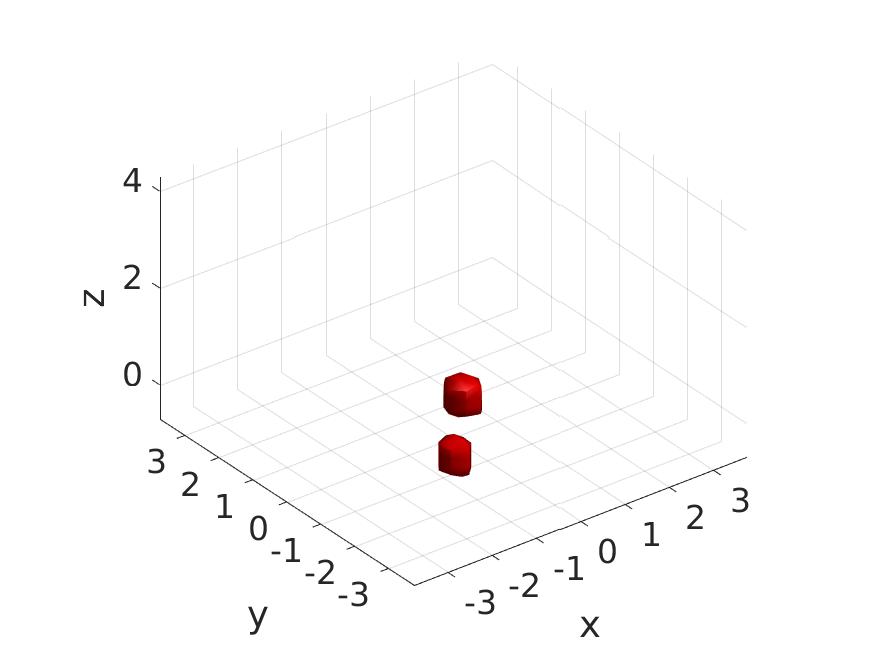}}\hfill
	\subfloat[\it Image computed from experimental data]{\includegraphics[width=0.3\linewidth]{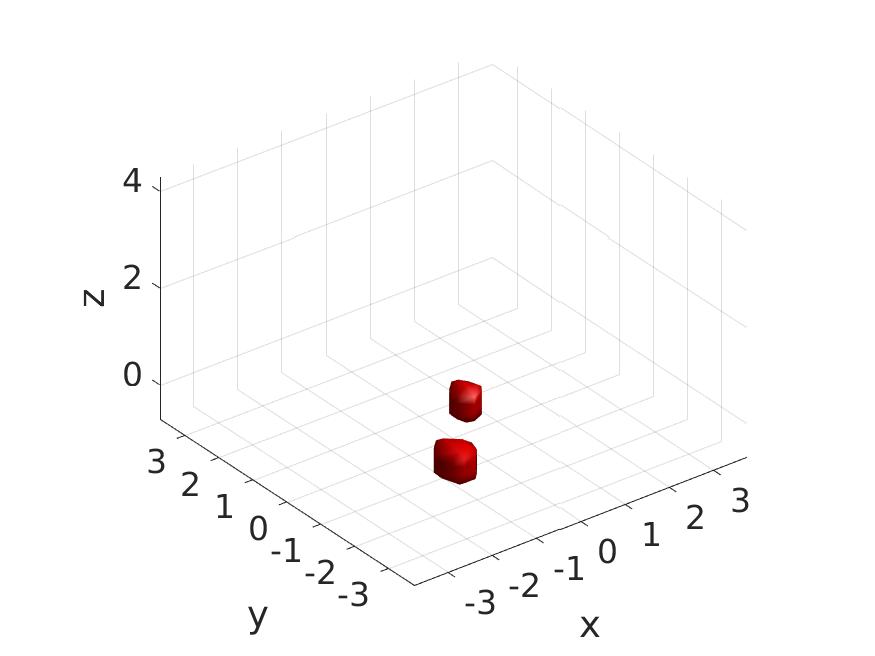}}
	\caption{\label{fig 9}\it Test 2. Correct and computed images of two microspheres involved in the experiment. One can observe that the locations of the unknown microspheres are imaged with a very good accuracy in both cases. A close similarity between (b) and (c) can also be observed.}
\end{figure}

\section{Summary}

\label{sec:8}

We have collected phaseless scattering data on six (6) wavelengths $\lambda $
listed in (\ref{6.1}) ranging from $\lambda =0.42\mu m$ to $\lambda
=0.671\mu m.$ To obtain the data at these wavelengths from the white light,
we have used narrow light filters. But only the data on $\lambda =0.525\mu m$
and $\lambda =0.580\mu m$ turned out to have a reasonable amount of noise.
So, we have linearly interpolated the measured data between points $k=%
\underline{k}=108.3$ and $k=\overline{k}=119.7,$ which mean corresponding
wave numbers. Only a single direction of the incident plane wave was used,
i.e. we have worked with the phaseless CIP with the
data resulting from a single measurement event. This is definitely more
challenging than the case of multiple measurements. We have measured the
intensity $\left\vert u\left( \mathbf{x},k\right) \right\vert ^{2}\mid _{P_{%
\text{meas}}}$of the full complex valued wave field on a square $P_{\text{%
meas}}$ located on a plane, which is orthogonal to the direction of the
propagation of the incident plane wave. Here $u\left( \mathbf{x},k\right) $
is the solution of the Helmholtz equation (\ref{eqn Helmholtz}) with the
radiation condition.

Since previous works on reconstruction procedures for phaseless CIPs \cite%
{KR:JIpp2015,KlibanovRomanov:ejmca2015,KlibanovRomanov:SIAMam2016,KlibanovRomanov:ip2016,KlibanovLocKejia:apnum2016,KlibanovLiemLoc:arxiv2017}%
of the first author with coauthors have discussed only phaseless CIPs with
measurements of $\left\vert u_{\text{sc}}\left( \mathbf{x},k\right)
\right\vert ^{2}$, where $u_{\text{sc}}\left( \mathbf{x},k\right) $ is the
scattered wave field, we have developed a new procedure to approximately
reconstruct the function $u\left( \mathbf{x},k\right) \mid _{P_{\text{meas}%
}} $ from its modulus measured on $P_{\text{meas}}.$

One of the key obstacles in this direction was the absence of a proper
analytical estimate of $\left\vert u_{\text{sc}}\left( \mathbf{x},k\right)
\right\vert ^{2}.$ While we have observed numerically that this term is
small indeed and can be dropped (Figures 2), it was unclear how to prove
this analytically. To obtain a proper estimate for this term, we have used
Theorems \ref{thm 1} and \ref{thm 4.1}. Both these theorems use ideas of the Riemannian
geometry and asymptotic analysis. While Theorem \ref{thm 1} was actually proved in 
\cite{KlibanovRomanov:ip2016}, Theorem \ref{thm 4.1} is completely new. As a result,
our upper estimate (\ref{4.5})\emph{\ }of $\left\vert u_{\text{sc}}\left( 
\mathbf{x},k\right) \right\vert ^{2}\mid _{P_{\text{meas}}}$is a reasonable
one for the given range of parameters. Furthermore, number-wise this
estimate is approximately the same as purely numerical estimates of Figures
\ref{fig 2}. Thus, the analytical estimate (\ref{4.5}) in combination with Figures \ref{fig 2}
ensure that the term $\left\vert u_{\text{sc}}\left( \mathbf{x},k\right)
\right\vert ^{2}\mid _{P_{\text{meas}}}$ is sufficiently small compared with
 1. Hence, this term can be dropped in (\ref{4.1414}). The
resulting formula (\ref{4.15}), in combination with Theorem \ref{thm 1}, enables us
to obtain the inversion formula \eqref{4.18}.

As soon as the function $u\left( \mathbf{x},k\right) \mid _{P_{\text{meas}}}$
is approximately reconstructed, one obtains the conventional phased CIP,
which, however, is also difficult to solve. Its numerical solution is the
second stage of our reconstruction procedure. To solve this problem
numerically, we have used the globally convergent numerical method which was
developed in \cite{KlibanovLiemLocHui:IPI2017} and then successfully tested
on microwave backscattering experimental data in \cite%
{AlekKlibanovLocLiemThanh:anm2017,LiemKlibanovLocAlekFiddyHui:jcp2017,nguyen2017:iip2017}%
.

In our numerical studies we have decided to verify the accuracy of our
inversions of experimental data via comparison of inversion results with
those of computationally simulated data. Thus, we have computationally
generated the data for exactly the same microspheres as we have used in
experiments. We have observed that computational results for the forward
problem have a very good similarity with experimental data, see Figures \ref{fig 4}, \ref{fig 5}, \ref{fig 8}. In addition, inversion results for both experimental data sets are
very similar with the those of computationally simulated data, see Figures \ref{fig 6}, \ref{fig 7}, \ref{fig 9}
and Tables \ref{table 1},\ref{table 2}. Finally, the reconstruction \ error of the refractive
index is between 5.2\% and 5.4\% in all cases, which is small.

Thus, we conclude that since results of the forward problem solution for
computationally simulated data are quite close to the experimental data and
also since our inversion provides quite accurate locations and refractive
indices of microspheres of interest for both types of data, then our
mathematical modeling of experimental data is quite accurate one, including
the drop of the term $\left\vert u_{\text{sc}}\left( \mathbf{x},k\right)
\right\vert ^{2}\mid _{P_{\text{meas}}}$ in (\ref{4.1414}).

\label{sec:9}

\begin{center}
\textbf{Acknowledgements}
\end{center}

The work of M. V. Klibanov and N. A. Koshev was supported by the US Army
Research Laboratory and US Army Research Office grant W911NF-15-1-0233 as
well as by the Office of Naval Research grant N00014-15-1-2330. The work of
L. H. Nguyen was partially supported by research funds FRG 111172 provided
by University of North Carolina at Charlotte.


\end{document}